\documentclass[12pt,draft]{article}

\usepackage[cp1251]{inputenc}
\usepackage[T2A]{fontenc}
\usepackage[english,ukrainian]{babel}
\usepackage{amsmath,amsfonts,amssymb}
\usepackage{geometry}
\usepackage{cite}
\newtheorem{theorem}{Теорема}
\newtheorem{lemma}{Лема}
\newtheorem{proposition}{Твердження}

\newtheorem{remark}{Зауваження}
\newtheorem{example}{Приклад}

\begin{document}

\begin{flushleft}

\textbf{Tetiana Kasirenko and Iryna Chepurukhina}\\
\small(Institute of Mathematics, National Academy of Sciences of Ukraine, Kyiv)

\medskip

\large\textbf{ELLIPTIC PROBLEMS IN THE SENSE OF LAWRUK\\ WITH BOUNDARY OPERATORS OF HIGHER ORDERS\\ IN REFINED SOBOLEV SCALE}\normalsize

\medskip\medskip

\textbf{Тетяна Касіренко і Ірина~Чепурухіна}\\
\small(Інститут математики НАН України, Київ)

\medskip

\large\textbf{ЕЛІПТИЧНІ ЗА ЛАВРУКОМ ЗАДАЧІ\\ З КРАЙОВИМИ ОПЕРАТОРАМИ ВИЩИХ ПОРЯДКІВ\\ В УТОЧНЕНІЙ СОБОЛЄВСЬКІЙ ШКАЛІ}\normalsize

\end{flushleft}

\medskip

\noindent In a refined Sobolev scale, we investigate an elliptic boundary--value problem with additional unknown functions in boundary conditions for which the maximum of orders of boundary operators is grater than or equal to the order of the elliptic equation. This scale consists of inner product H\"ormander spaces whose order of regularity is given by a real number and a function varying slowly at infinity in the sense of Karamata. We prove a theorem on the Fredholm property of a bounded operator corresponding to this problem in the refined Sobolev scale. For the generalized solutions to the problem, we establish a local a priory estimate and prove a theorem about their regularity in H\"ormander spaces. We find sufficient conditions under which given generalized derivatives of the solutions are continuous.

\medskip

\noindent В уточненій соболєвській шкалі досліджено еліптичну крайову задачу з додатковими невідомими функціями у крайових умовах, для яких максимум порядків крайових операторів більший за порядок еліптичного рівняння, або рівний йому. Ця шкала складається з гільбертових просторів Хермандера, для яких показником регулярності служать дійсне число і функція, повільно змінна на нескінченності за Караматою. Доведено теорему про нетеровість обмеженого оператора, відповідного цій задачі, в уточненій соболєвській шкалі. Для узагальнених розв'язків задачі встановлено локальну апріорну оцінку та доведено теорему про їх регулярність у просторах Хермандера. Знайдено достатні умови неперервності заданих узагальнених похідних розв'язків.

\bigskip

\section{Вступ}\label{1sec1} Ця стаття присвячена дослідженню еліптичних крайових задач з додатковими невідомими функціями у крайових умовах. Їх було введено Б.~Лавруком \cite{Lawruk63a, Lawruk63b, Lawruk65} у 1963 році. Вони природно виникають при переході від загальної (нерегулярної) еліптичної крайової до формально спряженої задачі відносно спеціальної формули Гріна. Клас еліптичних за Лавруком крайових задач є замкненим відносно такого переходу. Важливі приклади цих задач виникають у теорії пружності і гідродинаміці \cite{AslanyanVassilievLidskii81, Ciarlet90, NazarovPileckas93}.

Еліптичні за Лавруком крайові задачі досліджено у соболєвських просторах В.~О.~Козловим, В.~Г.~Маз'єю і Й.~Россманом \cite[розд.~3]{KozlovMazyaRossmann97} в основному для одного еліптичного рівняння та І.~Я.~Ройтберг \cite{RoitbergInna97, RoitbergInna98} для еліптичних систем мішаного порядку (див. також монографію \cite[розд.~2]{Roitberg99}). Було доведено теореми про нетеровість обмежених операторів, що відповідають цим задачам, і породжені ними ізоморфізми, та теореми про апріорні оцінки розв'язків задач і підвищення регулярності розв'язків.

Втім, соболєвська шкала є занадто грубою для низки важливих задач теорії рівнянь з частинними похідними \cite{Hermander65, Hermander86, MikhailetsMurach14, NicolaRodino10, Paneah00}. У цьому зв'язку Л.~Хермандер \cite{Hermander65, Hermander86} ще у 1963 році увів і дослідив широкий клас нормованих просторів, для яких показником регулярності служить не число (як у просторах Соболєва), а досить загальна вагова функція частотних змінних. Л.~Хермандер навів важливі застосування цих просторів до питань про розв'язність рівнянь з частинними похідними і регулярність їх розв'язків. Проте, довгий час простори Хермандера та їх різні версії не застосовували до крайових задач, оскільки не було виділено досить широких класів цих просторів, які б допускали коректне означення на гладких многовидах. Зауважимо, що серед просторів Хермандера найбільший інтерес з точки зору застосувань викликають саме гільбертові простори.

Недавно ситуація змінилася завдяки  роботам В.~А.~Михайлеця і О.~О.~Мурача \cite{MikhailetsMurach05UMJ5, MikhailetsMurach06UMJ3, MikhailetsMurach0607UMJ5, MikhailetsMurach06UMJ11, MikhailetsMurach06UMB4, MikhailetsMurach08UMJ4}, у яких побудовано загальну теорію розв'язності еліптичних крайових задач у класах гільбертових просторів Хермандера, що утворюють уточнену соболєвську шкалу. Показниками регулярності для цих просторів служить пара параметрів~--- числовий і функціональний. Останній є додатною функцією, повільно змінною на нескінченності за Й.~Караматою. Важливо, що уточнена соболєвська шкала отримується методом інтерполяції з функціональним параметром пар гільбертових просторів Соболєва. За його допомогою вдалося перенести основні результати “соболєвської” теорії еліптичних крайових задач на уточнену соболєвську шкалу (ці результати викладено у монографії \cite{MikhailetsMurach14} та огляді \cite{MikhailetsMurach12BJMA2}), а також на деякі більш широкі класи гільбертових просторів Хермандера \cite{Anop14Dop4, AnopMurach14UMJ7, AnopMurach14MFAT2, AnopKasirenko16MFAT4, KasirenkoMurach17UMJ11}.

У вказаних просторах Хермандера було досліджено і еліптичні за Лавруком крайові задачі \cite{Chepuruhina14Coll2, ChepurukhinaMurach15MFAT1, ChepurukhinaMurach15UMJ5, Chepuruhina15Coll2, Chepurukhina15DopNAN}. Втім, важливий випадок, коли максимум порядків крайових операторів більший за порядок еліптичного рівняння, або рівний йому, дотепер не вивчали у просторах Хермандера. Відмітимо, що еліптичні крайові задачі з крайовими операторами вищих порядків зустрічаються в акустиці, гідродинаміці, теорії випадкових процесів \cite{Venttsel59, Krasil'nikov61, VeshevKouzov77}.

Мета цієї роботи~--- дослідити в уточненій соболєвській шкалі характер розв'язності і властивості розв'язків еліптичних за Лавруком задач з крайовими операторами вищих порядків.

Робота складається з семи пунктів. Пункт~1 є вступом. У~п.~2 наведено постановку еліптичної за Лавруком крайової задачі. Там же розглянуто спеціальну формулу Гріна для цієї задачі у випадку крайових операторів вищих порядків і наведено формальну спряжену крайову задачу відносно зазначеної формули Гріна. Окрім того, розглянуто відповідні приклади.  У п.~3 наведено означення функціональних просторів Хермандера, які утворюють уточнену соболєвську шкалу. Пункт~4 містить основні результати роботи про властивості досліджуваної задачі у цій шкалі. Серед них~--- теореми про нетеровість оператора задачі у відповідних парах просторів Хермандера і породжені цією задачею ізоморфізми, теореми про локальну апріорну оцінку узагальнених розв'язків задачі та їх регулярність у просторах Хермандера. У п.~5, як застосування уточненої соболєвської шкали, отримано нові достатні умови неперервності узагальнених похідних розв'язків досліджуваної задачі, зокрема, умови класичності її узагальненого розв'язку. Усі результати роботи доведено у п.~6. Завершальний п.~7 містить висновки до роботи.

\section{Постановка задачі}\label{1sec2}

Нехай $\Omega$~--- довільна  обмежена  область у евклідовому просторі
$\mathbb{R}^{n}$, де $n\geq2$. Припустимо, що її межа $\Gamma:=\partial\Omega$ є нескінченно гладким замкненим (тобто компактним і без краю) многовидом вимірності $n-1$. При цьому вважаємо, що $C^{\infty}$-структура на $\Gamma$ породжена простором $\mathbb{R}^{n}$. Як звичайно, $\overline{\Omega}=\Omega\cup\Gamma$. Позначимо через $\nu(x)$ орт внутрішньої нормалі до межі $\Gamma$ у точці $x\in\Gamma$.

Довільно виберемо цілі числа $q\geq1$, $\varkappa\geq1$, $m_{1},\ldots,m_{q+\varkappa}$ і $r_{1},\ldots,r_{\varkappa}$.
Розглянемо в $\Omega$ таку лінійну крайову задачу:
\begin{gather}\label{1f1}
Au=f\quad\mbox{в}\quad\Omega,\\
B_{j}u+\sum_{k=1}^{\varkappa}C_{j,k}v_{k}=g_{j}\quad\mbox{на}\quad\Gamma,
\quad j=1,...,q+\varkappa.\label{1f2}
\end{gather}
Тут є невідомими функція $u$ в області $\Omega$ і $\varkappa$ функцій $v_{1},\ldots,v_{\varkappa}$ на $\Gamma$. В роботі всі функції та розподіли вважаємо комплекснозначними. У цій задачі
$$
A:=A(x,D):=\sum_{|\mu|\leq 2q}a_{\mu}(x)D^{\mu}
$$
є лінійним диференціальним оператором на $\overline{\Omega}$ парного порядку $2q\geq2$, кожне
$$
B_{j}(x,D)=\sum_{|\mu|\leq m_{j}}b_{j,\mu}(x)D^{\mu}
$$
є крайовим лінійним диференціальним оператором на $\Gamma$ порядку $\mathrm{ord}\,B_{j}\leq m_{j}$, а кожне $C_{j,k}=C_{j,k}(x,D_{\tau})$ є дотичним лінійним диференціальним оператором на $\Gamma$ порядку $\mathrm{ord}\,C_{j,k}\leq m_{j}+r_{k}$. При цьому, як звичайно, $B_{j}=0$, якщо $m_{j}<0$, і $C_{j,k}=0$, якщо $m_{j}+r_{k}<0$. Усі коефіцієнти цих диференціальних операторів є нескінченно гладкими функціями, заданими на $\overline{\Omega}$ і $\Gamma$ відповідно. Крайові задачі вигляду \eqref{1f1}, \eqref{1f2} уперше були розглянуті Б.~Лавруком \cite{Lawruk63a, Lawruk63b, Lawruk65}.

Тут і надалі використовуємо такі стандартні позначення:
$\mu:=(\mu_{1},\ldots,\mu_{n})$~--- мультиіндекс, $|\mu|:=\mu_{1}+\ldots+\mu_{n}$,
$D^{\mu}:=D_{1}^{\mu_{1}}\ldots D_{n}^{\mu_{n}}$, $D_{l}:=i\partial/\partial x_{l}$, де $l\in\{1,\ldots,n\}$, $i$~--- уявна одиниця, а $x=(x_1,\ldots,x_n)$~--- довільна точка простору $\mathbb{R}^{n}$. Окрім того, покладаємо $D_{\nu}:=i\partial/\partial\nu$ та $\xi^{\mu}:=\xi_{1}^{\mu_{1}}\ldots\xi_{n}^{\mu_{n}}$ для вектора $\xi=(\xi_{1},\ldots\xi_{n})\in\mathbb{C}^{n}$.

Припускаємо, що
$$
m:=\max\{m_{1},\ldots,m_{q+\varkappa}\}=
\max\{\mathrm{ord}\,B_{1},\ldots,\mathrm{ord}\,B_{q+\varkappa}\}.
$$
Окрім того, робимо природне припущення про те, що
\begin{equation}\label{nat-assump}
m\geq-r_{k}\quad\mbox{для кожного}\quad
k\in\{1,\ldots,\varkappa\}.
\end{equation}
(Якщо $m+r_{k}<0$ для деякого  $k$, то усі оператори $C_{1,k}$,..., $C_{q+\varkappa,k}$ дорівнюють нулю і тому шукана функція $v_{k}$ відсутня у крайових умовах \eqref{1f2}.)

Далі припускаємо, що крайова задача \eqref{1f1}, \eqref{1f2} є
еліптичною в області $\Omega$ за Лавруком. Наведемо відповідне означення (див., наприклад, \cite[п.~3.1.2]{KozlovMazyaRossmann97}).

Покладемо
$$
A^{\circ}(x,\xi):=\sum_{|\mu|=2q}a_{\mu}(x)\xi^{\mu}
\quad\mbox{для довільних}\quad x\in\overline{\Omega},\quad \xi\in\mathbb{C}^{n}.
$$
Вираз $A^{\circ}(x,\xi)$ є головним символом диференціального оператора $A(x,D)$. Аналогічно, для кожного номера $j\in\{1,\ldots,m+\varkappa\}$ покладемо
$$
B_{j}^{\circ}(x,\xi):=\sum_{|\mu|=m_{j}}b_{j,\mu}(x)\xi^{\mu}
\quad\mbox{для довільних}\quad x\in\Gamma,\quad \xi\in\mathbb{C}^{n};
$$
при цьому $B_{j}^{\circ}(x,\xi)\equiv0$, якщо $m_{j}<\mathrm{ord}\,B_{j}$. Якщо $m_{j}=\mathrm{ord}\,B_{j}$, то $B_{j}^{\circ}(x,\xi)$ є головним символом крайового диференціального оператора $B_{j}(x,D)$. Окрім того, для будь-яких номерів $j\in\{1,\ldots,m+\varkappa\}$ і $k\in\{1,\ldots,\varkappa\}$ позначимо через $C_{j,k}^{\circ}(x,\tau)$ головний символ дотичного диференціального оператора $C_{j,k}(x,D_{\tau})$, якщо $\mathrm{ord}\,C_{j,k}=m_{j}+r_{k}$. Для кожної точки $x\in\Gamma$ вираз $C^{\circ}_{j,k}(x,\tau)$ є однорідним поліномом порядку $m_{j}+r_{k}$ змінної $\tau$, де $\tau$--- довільний дотичний вектор до межі $\Gamma$ у точці~$x$. Якщо $\mathrm{ord}\,C_{j,k}<m_{j}+r_{k}$, то покладаємо $C^{\circ}_{j,k}(x,\tau):=0$.

Крайову задачу \eqref{1f1}, \eqref{1f2} називають еліптичною в області $\Omega$ за Лавруком, якщо виконуються такі три умови:
\begin{itemize}
\item [(i)] Диференціальний оператор $A(x,D)$ є еліптичним в кожній точці $x\in\overline{\Omega}$, тобто $A^{\circ}(x,\xi)\neq0$ для довільного вектора $\xi\in\mathbb{R}^{n}\setminus\{0\}$.
\item [(ii)] Диференціальний оператор $A(x,D)$ є правильно еліптичним в кожній точці $x\in\Gamma$, тобто для довільного вектора $\tau\neq0$, дотичного до $\Gamma$ у точці $x$, многочлен $A^{\circ}(x,\tau+\zeta\nu(x))$ комплексної змінної $\zeta$ має $q$ коренів з додатною уявною частиною і стільки ж коренів з від'ємною уявною частиною (підрахованих з урахуванням їх кратності).
\item [(iii)] Система крайових умов \eqref{1f2} накриває рівняння \eqref{1f1} у кожній точці $x\in\Gamma$, тобто для кожного вектора $\tau\neq0$, дотичного до $\Gamma$ у точці $x$, крайова задача
    \begin{equation}\label{elliptic-Lawruk-1}
    A^{\circ}(x,\tau+\nu(x)D_{t})\theta(t)=0\quad\mbox{при}\quad t>0,
    \end{equation}
    \begin{equation}\label{elliptic-Lawruk-2}
    \begin{gathered}
    B_{j}^{\circ}(x,\tau+\nu(x)D_{t})\theta(t)\big|_{t=0}+\\
    +\sum_{k=1}^{\varkappa}C_{j,k}^{\circ}(x,\tau)\lambda_{k}=0,\quad j=1,...,q+\varkappa,
    \end{gathered}
    \end{equation}
    має лише тривіальний (нульовий) розв'язок. Ця задача розглядається відносно невідомої функції $\theta\in C^{\infty}([0,\infty))$, яка задовольняє умову $\theta(t)\rightarrow0$ при $t\rightarrow\infty$, і невідомих комплексних чисел $\lambda_{1},\ldots,\lambda_{\varkappa}$. Тут $A^{\circ}(x,\tau+\nu(x)D_{t})$  і $B_{j}^{\circ}(x,\tau+\nu(x)D_{t})$ є диференціальними операторами відносно $D_{t}:=i\partial/\partial t$, які отримуємо, поклавши $\zeta:=D_{t}$ у відповідно многочленах $A^{\circ}(x,\tau+\zeta\nu(x))$ і $B_{j}^{\circ}(x,\tau+\zeta\nu(x))$ змінної~$\zeta$.
\end{itemize}

Зауважимо \cite[с.~166]{FunctionalAnalysis72}, що умова (ii) є наслідком умови (i) у випадку, коли $n\geq3$, а також у випадку, коли $n=2$ і усі старші коефіцієнти оператора $A(x,D)$ є дійсними функціями.

\begin{example}\label{1ex1} \rm
Нехай $n=2$, $q=1$ і $\varkappa=1$. Припустимо, що диференціальний оператор $A(x,D)$ другого порядку задовольняє умови (i) та (ii). Розглянемо крайову задачу, яка складається з диференціального рівняння \eqref{1f1} і пари крайових умов
\begin{gather*}
\partial_{\nu}^{p}u+v=g_{1}\quad\mbox{на}\quad\Gamma,\\
\partial_{\nu}^{p+1}u+\partial_{\Gamma}v=g_{2}\quad\mbox{на}\quad\Gamma,
\end{gather*}
де довільно вибрано ціле $p\geq0$. Тут  $\partial_{\nu}:=\partial/\partial\nu$~--- похідна вздовж орта $\nu$,
а $\partial_{\Gamma}$~--- похідна вздовж кривої~$\Gamma$ у додатному напряму. Ці крайові умови набирають вигляду  \eqref{1f2}, де $B_{1}=\partial_{\nu}^{p}$, $B_{2}=\partial_{\nu}^{p+1}$ і $C_{1,1}=1$, $C_{2,1}=\partial_{\Gamma}$, а $m_{1}=p$, $m_{2}=p+1$ і $r_{1}=-p$.

Переконаємося, що розглянута крайова задача є еліптичною за Лавруком в області~$\Omega$. Треба перевірити лише, що вона задовольняє умову~(iii). Виберемо довільно точку $x\in\Gamma$ і вектор $\tau\neq0$, дотичний до $\Gamma$ у цій точці. Загальний розв'язок диференціального рівняння \eqref{elliptic-Lawruk-1}, який задовольняє умову $\theta(t)\rightarrow0$ при $t\rightarrow\infty$, записується у вигляді
\begin{equation}\label{ex3.3-theta}
\theta(t)=c\exp(-i\zeta_{-}(x,\tau)t),
\end{equation}
де $c$ є довільне комплексне число, а $\zeta_{-}(x,\tau)$ є $\zeta$-корінь многочлена $A^{\circ}(x,\tau+\zeta\nu(x))$ такий, що $\mathrm{Im}\,\zeta_{-}(x,\tau)<0$. Для розглянутої крайової задачі умови \eqref{elliptic-Lawruk-2} набирають вигляду
\begin{gather}\label{ex3.3-1a}
\theta^{(p)}(0)+\lambda_{1}=0, \\
\theta^{(p+1)}(0)\mp i|\tau|\lambda_{1}=0. \label{ex3.3-1b}
\end{gather}

Справді, оскільки
\begin{equation*}
B_{j}(x,D)=\partial_{\nu}^{p+j-1}=(-i)^{p+j-1}D^{p+j-1}_{\nu}
\end{equation*}
для кожного номера $j\in\{1,2\}$, то в умові (iii)
\begin{equation*}
B_{j}^{\circ}(x,\tau+\zeta\nu(x))=(-i\zeta)^{p+j-1}.
\end{equation*}
Отже,
\begin{equation}\label{ex3.3-B2}
B_{j}^{\circ}(x,\tau+\nu(x)D_{t})\theta(t)=
(-iD_{t})^{p+j-1}\theta(t)=\theta^{(p+j-1)}(t).
\end{equation}
Окрім того,
\begin{equation*}
C_{2,1}(x,D_{\Gamma})=\partial_{\Gamma}=
\frac{\pm\tau_{1}}{|\tau|}\frac{\partial}{\partial x_{1}}+
\frac{\pm\tau_{2}}{|\tau|}\frac{\partial}{\partial x_{2}}=
\frac{\mp i}{|\tau|}(\tau_{1}D_{1}+\tau_{2}D_{2}),
\end{equation*}
де записуємо $\tau=(\tau_{1},\tau_{2})$ і вибираємо верхній (відповідно нижній) знак, якщо дотичний вектор $\tau$ направлений у бік додатного (від'ємного) обходу кривої $\Gamma$. Тому
\begin{equation}\label{ex3.3-C2}
C_{2,1}^{\circ}(x,\tau)=
\frac{\mp i}{|\tau|}(\tau_{1}^{2}+\tau_{2}^{2})
=\mp i|\tau|.
\end{equation}
З огляду на \eqref{ex3.3-B2} і \eqref{ex3.3-C2} робимо висновок, що крайова умова \eqref{elliptic-Lawruk-2} при $j=1$ набирає вигляду \eqref{ex3.3-1a}, а при $j=2$~--- вигляду \eqref{ex3.3-1b}.

Підставивши \eqref{ex3.3-theta} в умови \eqref{ex3.3-1a} і \eqref{ex3.3-1b}, запишемо
\begin{gather*}
(-i\zeta_{-}(x,\tau))^{p}c+\lambda_{1}=0, \\
(-i\zeta_{-}(x,\tau))^{p+1}c\mp i|\tau|\lambda_{1}=0.
\end{gather*}
Це~--- система лінійних однорідних алгебраїчних рівнянь відносно невідомих $c,\lambda_{1}\in\mathbb{C}$. Її визначник
\begin{equation*}
\left|
  \begin{array}{cc}
    (-i\zeta_{-}(x,\tau))^{p} & 1 \\
    (-i\zeta_{-}(x,\tau))^{p+1} & \mp i|\tau| 
  \end{array}
\right|=-i(-i\zeta_{-}(x,\tau))^{p}(\pm |\tau|-\zeta_{-}(x,\tau))
\end{equation*}
відмінний від нуля. Тому ця система має лише тривіальний розв'язок, тобто умова (iii) виконується. Отже, розглянута у цьому прикладі крайова задача є еліптичною за Лавруком в області~$\Omega$.
\end{example}

\begin{example} \rm
Нехай $n\geq2$, $q=2$ і $\varkappa=2$. Припустимо, що диференціальний оператор $A(x,D)$ четвертого порядку задовольняє умови (i) та (ii). Розглянемо крайову задачу, яка складається з диференціального рівняння \eqref{1f1} і таких чотирьох крайових умов на~$\Gamma$:
\begin{gather*}
\partial_{\nu}^{p}u+v_{1}=g_{1},\\
\partial_{\nu}^{p+1}u+v_{2}=g_{2},\\
\partial_{\nu}^{p+2}u+\Delta_{\Gamma}v_{1}=g_{3},\\
\partial_{\nu}^{p+3}u+\Delta_{\Gamma}v_{2}=g_{4}.
\end{gather*}
Тут ціле $p\geq0$ вибрано довільно, а $\Delta_{\Gamma}$ позначає, як звичайно, оператор Бельтрамі-Лапласа на $\Gamma$, при цьому на $\Gamma$ введено ріманову метрику, індуковану простором $\mathbb{R}^{n}$. Ці крайові умови набирають вигляду  \eqref{1f2}, де $B_{j}=\partial_{\nu}^{p+j-1}$ для кожного $j\in\{1,2,3,4\}$, а $C_{1,1}=1$, $C_{1,2}=0$, $C_{2,1}=0$, $C_{2,2}=1$, $C_{3,1}=\Delta_{\Gamma}$, $C_{3,2}=0$, $C_{4,1}=0$ і $C_{4,2}=\Delta_{\Gamma}$. Тут $m_{j}=p+j-1$ для кожного $j\in\{1,2,3,4\}$, а $r_{1}=-p$ і $r_{2}=-p-1$.

Додатково припустимо, що для довільної точки $x\in\Gamma$ і вектора $\tau\neq0$, дотичного до $\Gamma$ у цій точці, многочлен $A^{\circ}(x,\tau+\zeta\nu(x))$ не має кратних $\zeta$-коренів. Це припущення виконується, наприклад, якщо $A=\partial_{\nu}^{4}+\Delta_{\Gamma}^{2}$ на~$\Gamma$, бо тоді $A^{\circ}(x,\tau+\zeta\nu(x))=\zeta^{4}+\gamma^{2}(x,\tau)$. Тут $\gamma(x,\tau)$ позначає головний символ оператора $\Delta_{\Gamma}$; цей символ задовольняє умову $\gamma(x,\tau)<0$ при $\tau\neq0$.

Перевіримо, що розглянута крайова задача є еліптичною за Лавруком в області $\Omega$. Треба показати лише, що вона задовольняє умову~(iii). Виберемо довільно точку $x\in\Gamma$ і вектор $\tau\neq0$, дотичний до $\Gamma$ у цій точці. Для цієї задачі система крайових умов~\eqref{elliptic-Lawruk-2} набирає вигляду
\begin{gather*}
\theta^{(p)}(0)+\lambda_{1}=0,\\
\theta^{(p+1)}(0)+\lambda_{2}=0,\\
\theta^{(p+2)}(0)+\gamma(x,\tau)\lambda_{1}=0,\\
\theta^{(p+3)}(0)+\gamma(x,\tau)\lambda_{2}=0.
\end{gather*}
Тут $\theta(t)$~--- загальний розв'язок диференціального рівняння \eqref{1f1}, який задовольняє умову $\theta(t)\rightarrow0$ при $t\rightarrow\infty$, а $\lambda_{1}, \lambda_{2}$~--- довільні комплексні числа. Ця система еквівалентна такій:
\begin{gather}
\lambda_{1}=-\theta^{(p)}(0), \notag \\
\quad \lambda_{2}=-\theta^{(p+1)}(0),\notag \\
\theta^{(p+2)}(0)-\gamma(x,\tau)\theta^{(p)}(0)=0,\label{ex2-a}\\
\theta^{(p+3)}(0)-\gamma(x,\tau)\theta^{(p+1)}(0)=0.\label{ex2-b}
\end{gather}

Згідно зроблених припущень щодо оператора $A$ розв'язок $\theta(t)$ запишемо у вигляді
\begin{equation*}
\theta(t)=c_{1}\exp(-i\zeta_{1}(x,\tau)t)+c_{2}\exp(-i\zeta_{2}(x,\tau)t),
\end{equation*}
де $\zeta_{1}(x,\tau)$ і $\zeta_{2}(x,\tau)$~--- різні $\zeta$-корені многочлена $A^{(0)}(x,\tau+\zeta\nu(x))$, які мають від'ємну уявну частину, а $c_{1}$ і $c_{2}$~--- довільні комплексні числа. Підставивши цей розв'язок в умови \eqref{ex2-a} і \eqref{ex2-b}, отримаємо систему однорідних лінійних алгебраїчних рівнянь
\begin{gather*}
(\zeta_{1}(x,\tau))^{p}\gamma_{1}(x,\tau)c_{1}+
(\zeta_{2}(x,\tau))^{p}\gamma_{2}(x,\tau)c_{2}=0,\\
(\zeta_{1}(x,\tau))^{p+1}\gamma_{1}(x,\tau)c_{1}+
(\zeta_{2}(x,\tau))^{p+1}\gamma_{2}(x,\tau)c_{2}=0,
\end{gather*}
відносно невідомих чисел $c_{1}$ і $c_{2}$, де позначено
\begin{equation*}
\gamma_{j}(x,\tau):=-\zeta_{j}^{2}(x,\tau)-\gamma(x,\tau)\quad
\mbox{для кожного}\quad j\in\{1,2\}.
\end{equation*}
Отже, для розглянутої у цьому прикладі крайової задачі умова (iii) рівносильна тому, що ця система має лише тривіальний розв'язок, тобто її головний визначник
\begin{equation*}
(\zeta_{1}(x,\tau))^{p}(\zeta_{2}(x,\tau))^{p}\gamma_{1}(x,\tau)
\gamma_{2}(x,\tau)(\zeta_{2}(x,\tau)-\zeta_{1}(x,\tau))\neq0.
\end{equation*}
Оскільки корені $\zeta_{1}(x,\tau)$ і $\zeta_{2}(x,\tau)$ різні та ненульові, то остання умова еквівалентна такій
\begin{equation}\label{ex3.5-f5}
\gamma_{1}(x,\tau)\gamma_{2}(x,\tau)\neq0.
\end{equation}

Покажемо, що ця нерівність виконується.  Для кожного номера $j\in\{1,2\}$ запишемо $\zeta_{j}(x,\tau)=\alpha_{j}+i\beta_{j}$, де $\alpha_{j}\in\mathbb{R}$ та $\beta_{j}<0$. Тоді
\begin{equation*}
\gamma_{j}(x,\tau)=-(\alpha_{j}+i\beta_{j})^{2}-\gamma(x,\tau)=
\beta_{j}^{2}-\alpha_{j}^{2}-\gamma(x,\tau)-2i\alpha_{j}\beta_{j};
\end{equation*}
тут, нагадаємо, $\gamma(x,\tau)<0$. Таким чином, якщо $\alpha_{j}=0$, то $\gamma_{j}(x,\tau)>0$, а якщо $\alpha_{j}\neq0$, то $\gamma_{j}(x,\tau)\notin\mathbb{R}$. Отже, умова \eqref{ex3.5-f5} виконується і тому розглянута крайова задача є еліптичною за Лавруком в області $\Omega$.
\end{example}

Далі припускаємо, що $m\geq2q$, тобто принаймні один крайовий оператор $B_{j}$ має порядок $\mathrm{ord}\,B_{j}\geq2q$.

Пов'яжемо із задачею \eqref{1f1}, \eqref{1f2} лінійне відображення
\begin{equation}\label{1f3}
\begin{gathered}
\Lambda:(u,v_{1},...,v_{\varkappa})\rightarrow\\ \rightarrow
\biggl(Au,\,B_{1}u+\sum_{k=1}^{\varkappa}C_{1,k}v_{k},...,
B_{q+\varkappa}u+\sum_{k=1}^{\varkappa}C_{q+\varkappa,k}v_{k}\biggr),\\
\mbox{де}\quad u\in C^{\infty}(\overline{\Omega}),\quad
v_{1},\ldots,v_{\varkappa}\in C^{\infty}(\Gamma).
\end{gathered}
\end{equation}
Ми досліджуємо властивості продовження за неперервністю цього відображення у підходящих парах гільбертових просторів Хермандера, які утворюють уточнену соболєвську шкалу.

Для опису області значень цього продовження нам потрібна така спеціальна формула Гріна \cite[формула (4.1.10)]{KozlovMazyaRossmann97}:
\begin{gather*}
(Au,\omega)_{\Omega}+\sum_{j=1}^{m-2q+1}(D_{\nu}^{j-1}Au,w_{j})_{\Gamma}+
\sum_{j=1}^{q+\varkappa}\biggl(B_{j}u+
\sum_{k=1}^{\varkappa}C_{j,k}v_{k},h_{j}\biggr)_{\Gamma}=\\
=(u,A^{+}\omega)_{\Omega}+\sum_{k=1}^{m+1}\biggl(D_{\nu}^{k-1}u,K_{k}\omega+
\sum_{j=1}^{m-2q+1}R_{j,k}^{+}w_{j}+
\sum_{j=1}^{q+\varkappa}Q_{j,k}^{+}h_{j}\biggr)_{\Gamma}+\sum_{k=1}^{\varkappa}\biggl(v_{k},
\sum_{j=1}^{q+\varkappa}C_{j,k}^{+}h_{j}\biggr)_{\Gamma},
\end{gather*}
для довільних функцій $u,\omega\in C^{\infty}(\overline{\Omega})$ і
$$
v_{1},\ldots,v_{\varkappa},w_{1},\ldots,w_{m-2q+1},
h_{1},\ldots,h_{q+\varkappa}\in C^{\infty}(\Gamma).
$$
Тут і далі через $(\cdot,\cdot)_{\Omega}$ і $(\cdot,\cdot)_{\Gamma}$ позначено відповідно скалярні добутки у гільбертових просторах $L_{2}(\Omega)$ і $L_{2}(\Gamma)$ функцій квадратично інтегровних на $\Omega$ і $\Gamma$ відносно мір Лебега, а також продовженням за неперервністю цих скалярних добутків. Як звичайно, $A^{+}$~--- диференціальний оператор, формально спряжений до $A$, тобто
$$
(A^{+}\omega)(x):=
\sum_{|\mu|\leq2q}D^{\mu}(\overline{a_{\mu}(x)}\omega(x)).
$$
Окрім того, усі $C_{j,k}^{+}$, $R_{j,k}^{+}$ і $Q_{j,k}^{+}$ є дотичними диференціальними операторами, формально спряженими  відповідно до $C_{j,k}$, $R_{j,k}$ і $Q_{j,k}$ відносно $(\cdot,\cdot)_{\Gamma}$, а дотичні лінійні диференціальні оператори $R_{j,k}:=R_{j,k}(x,D_{\tau})$ і $Q_{j,k}:=Q_{j,k}(x,D_{\tau})$ узяті із зображення крайових диференціальних операторів $D_{\nu}^{j-1}A$ і $B_{j}$ у вигляді
\begin{gather*}
D_{\nu}^{j-1}A(x,D)=\sum_{k=1}^{m+1}R_{j,k}(x,D_{\tau})D_{\nu}^{k-1},\quad
j=1,\ldots,m-2q+1,\\
B_{j}(x,D)=\sum_{k=1}^{m+1}Q_{j,k}(x,D_{\tau})D_{\nu}^{k-1},\quad
j=1,\ldots,q+\varkappa.
\end{gather*}
Відмітимо, що $\mathrm{ord}\,R_{j,k}\leq 2q+j-k$ і $\mathrm{ord}\,Q_{j,k}\leq m_{j}-k+1$; при цьому, звісно, $R_{j,k}=0$ при $k\geq2q+j+1$ і $Q_{j,k}=0$ при $k\geq m_{j}+2$. Нарешті, кожне $K_{k}:=K_{k}(x,D)$~--- деякий крайовий лінійний диференціальний оператор на $\Gamma$ порядку $\mathrm{ord}\,K_{k}\leq2q-k$ з коефіцієнтами класу $C^{\infty}(\Gamma)$; при цьому, якщо $k\geq2q+1$, то $K_{k}=0$.

Спеціальна формула Гріна приводить до такої крайової
задачі в області $\Omega$:
\begin{equation}\label{1f4}
A^{+}\omega=\eta\quad\mbox{в}\quad\Omega,
\end{equation}
\begin{equation}\label{1f5}
\begin{gathered}
K_{k}\omega+\sum_{j=1}^{m-2q+1}R_{j,k}^{+}w_{j}+
\sum_{j=1}^{q+\varkappa}Q_{j,k}^{+}h_{j}=\psi_{k}\quad
\mbox{на}\quad\Gamma,\\\quad k=1,...,m+1,
\end{gathered}
\end{equation}
\begin{equation}\label{1f6}
\sum_{j=1}^{q+\varkappa}C_{j,k}^{+}h_{j}=\psi_{m+1+k}
\quad\mbox{на}\quad\Gamma,\quad k=1,...,\varkappa.
\end{equation}
Ця задача містить, окрім невідомої функції $\omega$ в області $\Omega$, ще $m-q+1+\varkappa$ додаткових невідомих функцій $w_{1},\ldots,w_{m-2q+1}$ і $h_{1},\ldots,h_{q+\varkappa}$ на межі $\Gamma$. Задачу \eqref{1f4}\,--\,\eqref{1f6} називають формально спряженою до задачі \eqref{1f1}, \eqref{1f2} відносно розглянутої спеціальної формули Гріна.  Відомо \cite[теорема~4.1.1]{KozlovMazyaRossmann97}, що еліптичність за Лавруком задачі \eqref{1f1}, \eqref{1f2} рівносильна еліптичності за Лавруком формально спряженої задачі \eqref{1f4}\,--\,\eqref{1f6}.

\begin{example}\label{1ex1-formally-adj} \rm
Запишемо спеціальну формулу Гріна для крайової задачі
\begin{equation}\label{51f1ex2}
\begin{gathered}
\Delta u=f\quad\mbox{в}\quad\Omega,\\
\partial_{\nu}u+v=g_{1},\quad
\partial_{\nu}^{2}u+\partial_{\Gamma}v=g_{2}\quad\mbox{на}\quad\Gamma,
\end{gathered}
\end{equation}
заданої в крузі
$$
\Omega:=\{(x_1,x_2)\in \mathbb{R}^{2}: x_{1}^2+x_{2}^2<1\}.
$$
Ця задача є еліптичною за Лавруком в $\Omega$, оскільки вона є окремим випадком еліптичної крайової задачі, розглянутої у прикладі~\ref{1ex1}.
Зауважимо, що $\Delta u=\partial_{\nu}^{2}u-\partial_{\nu}u+\partial_{\Gamma}^{2}u$ на~$\Gamma$; тут $\partial_{\nu}=-\partial/\partial\varrho$ і $\partial_{\Gamma}=\partial/\partial\varphi$, а $(\varrho,\varphi)$~--- полярні координати. Застосувавши другу класичну формулу Гріна для оператора Лапласа, отримаємо такі рівності:
\begin{gather*}
(\Delta u,\omega)_{\Omega}+(\Delta u,w)_{\Gamma}+(\partial_{\nu}u+v,h_{1})_{\Gamma}
+(\partial_{\nu}^{2}u+\partial_{\Gamma}v,h_{2})_{\Gamma}=\\
=(u,\Delta\omega)_{\Omega}-(\partial_{\nu}u,\omega)_{\Gamma}
+(u,\partial_{\nu}\omega)_{\Gamma}+(\partial_{\nu}^{2}u-
\partial_{\nu}u+\partial_{\Gamma}^{2}u,w)_{\Gamma}+\\
+(\partial_{\nu}u+v,h_{1})_{\Gamma}
+(\partial_{\nu}^{2}u+\partial_{\Gamma}v,h_{2})_{\Gamma}=\\
=(u,\Delta\omega)_{\Omega}+
(u,\partial_{\nu}\omega+\partial_{\Gamma}^{2}w)_{\Gamma}+
(\partial_{\nu}u,-\omega-w+h_{1})_{\Gamma}+\\
+(\partial_{\nu}^{2}u,w+h_{2})_{\Gamma}+
(v,h_{1}-\partial_{\Gamma}h_{2})_{\Gamma}
\end{gather*}
для довільних функцій $u,\omega\in C^{\infty}(\overline{\Omega})$ і $v,w,h_{1},h_{2}\in C^{\infty}(\Gamma)$. Отже, спеціальна формула Гріна для крайової задачі \eqref{51f1ex2} набирає вигляду
\begin{gather*}
(\Delta u,\omega)_{\Omega}+
(\Delta u,w)_{\Gamma}+(\partial_{\nu}u+v,h_{1})_{\Gamma}
+(\partial_{\nu}^{2}u+\partial_{\Gamma}v,h_{2})_{\Gamma}=\\
=(u,\Delta\omega)_{\Omega}+
(u,\partial_{\nu}\omega+\partial_{\Gamma}^{2}w)_{\Gamma}+
(D_{\nu}u,-i\omega-iw+ih_{1})_{\Gamma}+\\
+(D_{\nu}^{2}u,-w-h_{2})_{\Gamma}+(v,h_{1}-\partial_{\Gamma}h_{2})_{\Gamma}.
\end{gather*}
Тому крайова задача
\begin{gather*}
\Delta\omega=\eta\quad\mbox{в}\quad\Omega,\\
\partial_{\nu}\omega+\partial_{\Gamma}^{2}w=\psi_{1},\quad
-i\omega-iw+ih_{1}=\psi_{2}\quad\mbox{на}\quad\Gamma,\\
-w-h_{2}=\psi_{3},\quad
h_{1}-\partial_{\Gamma}h_{2}=\psi_{4}\quad\mbox{на}\quad\Gamma
\end{gather*}
є формально спряженою до задачі \eqref{51f1ex2} відносно цієї формули Гріна. Отримана формально спряжена задача містить три додаткові невідомі функції $w$, $h_{1}$ і $h_{2}$ на $\Gamma$.
\end{example}

\section{Уточнена соболєвська шкала}\label{1sec3}

Розглянемо гільбертові простори Хермандера \cite[п.~2.2]{Hermander65}, в яких будемо досліджувати крайову задачу \eqref{1f1}, \eqref{1f2}. Вони утворюють уточнені соболєвські шкали $\{H^{s,\varphi}(G):s\in\mathbb{R},\varphi\in\mathcal{M}\}$, де $G\in\{\Omega,\Gamma\}$, введені та досліджені В.~А.~Михайлецем і О.~О.~Мурачем \cite{MikhailetsMurach05UMJ5, MikhailetsMurach06UMJ3} (див. також їх монографію \cite{MikhailetsMurach14} та огляд
\cite{MikhailetsMurach12BJMA2}). Тут і надалі
$\mathcal {M}$~--- множина всіх вимірних за Борелем функцій $\varphi: [1,\infty)\rightarrow(0,\infty)$, які обмежені і відокремлені від нуля на кожному компакті і повільно змінюються на нескінченності за Караматою~\cite{Karamata30a}. Остання властивість значить, що $\varphi(\lambda t)/\varphi(t)\rightarrow 1$ при $t\rightarrow\infty$ для кожного $\lambda>0$. Повільно змінні функції добре вивчені і мають різноманітні застосування (див., наприклад, монографії \cite{Seneta85, BinghamGoldieTeugels89}).

Характерним прикладом функції класу $\mathcal {M}$ служить неперервна функція $\varphi: [1,\infty)\rightarrow(0,\infty)$ така, що
$$
\varphi(t)=(\log t)^{r_{1}}(\log\log
t)^{r_{2}}\ldots(\underbrace{\log\ldots\log}_{k\;\mbox{\small разів}}
t)^{r_{k}}\quad\mbox{при}\quad t\gg1;
$$
тут числа $k\in\mathbb{N}$ і $r_{1},\ldots,r_{k}\in\mathbb{R}$ вибрано довільно.

Нехай $s\in\mathbb{R}$ і $\varphi\in\mathcal {M}$. Означимо простір
$H^{s,\varphi}(\cdot)$ спочатку на $\mathbb {R}^{n}$, а потім на $\Omega$ і $\Gamma$. Будемо дотримуватися монографії \cite[пп.~1.3, 2.1, 3.2]{MikhailetsMurach14}.

За означенням, комплексний лінійний простір $H^{s,\varphi}(\mathbb{R}^{n})$, де $n\in\mathbb{N}$, складається з усіх повільно зростаючих розподілів $w$ на $\mathbb{R}^{n}$ таких, що їх перетворення Фур'є $\widehat{w}$ є функцією, локально інтегровною на $\mathbb {R}^{n}$ за Лебегом, і виконується умова
$$
\int\limits_{\mathbb{R}^{n}}
\langle\xi\rangle^{2s}\varphi^{2}(\langle\xi\rangle)|
\widehat{w}(\xi)|^{2}d\xi<\infty.
$$
Тут $\langle\xi\rangle:=(1+|\xi|^{2})^{1/2}$ є згладженим модулем вектора $\xi\in\mathbb{R}^{n}$. У~просторі $H^{s,\varphi}(\mathbb {R}^{n})$ введений скалярний добуток розподілів $w_{1}$ і $w_{2}$ за формулою
$$
(w_{1},w_{2})_{H^{s,\varphi}(\mathbb {R}^{n})}:= \int\limits_{\mathbb{R}^{n}}
\langle\xi\rangle^{2s}\varphi^{2}(\langle\xi\rangle)
\widehat{w_{1}}(\xi)\overline{\widehat{w_{2}}(\xi)}d\xi.
$$
Він породжує норму
$$
\|w\|_{H^{s,\varphi}(\mathbb{R}^{n})}:=
(w,w)_{H^{s,\varphi}(\mathbb{R}^{n})}^{1/2}.
$$

Простір $H^{s,\varphi}(\mathbb{R}^{n})$ є ізотропним гільбертовим випадком простору $\mathcal{B}_{p,\mu}$, введеного і дослідженого Л.~Хермандером \cite[п.~2.2]{Hermander65} (див. також його монографію \cite[п.~10.1]{Hermander86}). А саме, $H^{s,\varphi}(\mathbb{R}^{n})=\mathcal{B}_{2,\mu}$, якщо  $\mu(\xi)=\langle\xi\rangle^{s}\varphi(\langle\xi\rangle)$ для довільного $\xi\in\mathbb{R}^{n}$. Зауважимо, що у гільбертовому випадку $p=2$ простір $\mathcal{B}_{p,\mu}$ і його версії для евклідових областей ввели і дослідили також Л.~Р.~Волевич і Б.~П.~Панеях \cite[\S~2]{VolevichPaneah65}.

У важливому окремому випадку $\varphi(\cdot)\equiv1$ простір $H^{s,\varphi}(\mathbb {R}^{n})$ стає гільбертовим простором Соболєва $H^{s}(\mathbb {R}^{n})$ порядку $s$. У загальній ситуації виконуються неперервні та щільні вкладення
\begin{equation}\label{1f7}
\begin{gathered}
H^{s+\varepsilon}(\mathbb{R}^{n})\hookrightarrow H^{s,\varphi}(\mathbb{R}^{n})\hookrightarrow H^{s-\varepsilon}(\mathbb{R}^{n})\\
\mbox{для довільного}\quad\varepsilon>0.
\end{gathered}
\end{equation}
З них видно, що у класі функціональних просторів
\begin{equation*}
\bigl\{H^{s,\varphi}(\mathbb{R}^{n}):
s\in\mathbb{R},\varphi\in\mathcal {M}\bigr\}
\end{equation*}
числовий параметр $s$ задає основну (степеневу) регулярність, а функціональний параметр $\varphi$ задає додаткову регулярність,  підпорядковану основній. У залежності від того, чи $\varphi(t)\rightarrow\infty$ або $\varphi(t)\rightarrow 0$ при $t\rightarrow\infty$, параметр $\varphi$ задає додаткову додатну або від'ємну регулярність. Іншими словами, параметр $\varphi$ уточнює основну $s$-гладкість. Тому цей клас природно називати уточненою соболєвською шкалою на $\mathbb {R}^{n}$.

Її аналоги для евклідової області $\Omega$ і замкненого
компактного многовиду $\Gamma$ вводяться у стандартний спосіб.
Наведемо відповідні означення.

Комплексний лінійний простір $H^{s,\varphi}(\Omega)$ складається зі звужень в
область $\Omega$ усіх розподілів $w\in H^{s,\varphi}(\mathbb{R}^{n})$.
Норма у ньому означена за формулою
$$
\|u\|_{H^{s,\varphi}(\Omega)}:=
\inf\bigl\{\,\|w\|_{H^{s,\varphi}(\mathbb{R}^{n})}:\,
w\in
H^{s,\varphi}(\mathbb{R}^{n}),\;w=u\;\,\mbox{в}\;\,\Omega\,\bigr\},
$$
де $u\in H^{s,\varphi}(\Omega)$. Цей простір гільбертів і сепарабельний відносно вказаної норми. Множина $C^{\infty}(\overline{\Omega})$ щільна у ньому.

Комплексний лінійний простір $H^{s,\varphi}(\Gamma)$ складається з усіх розподілів на $\Gamma$, які в локальних координатах дають елементи простору $H^{s,\varphi}(\mathbb{R}^{n-1})$. Дамо докладне означення. Нехай довільним чином вибрано скінченний атлас із $C^{\infty}$-структури на многовиді $\Gamma$, утворений локальними картами $\pi_j:
\mathbb{R}^{n-1}\leftrightarrow\Gamma_{j}$, де $j=1,\ldots,\lambda$.
Тут відкриті множини $\Gamma_{1},\ldots,\Gamma_{\lambda}$ складають скінченне покриття многовиду $\Gamma$. Нехай, окрім того, вибрані функції $\chi_j\in C^{\infty}(\Gamma)$, де $j=1,\ldots,\lambda$, які утворюють розбиття одиниці на $\Gamma$, що задовольняє умову $\mathrm{supp}\,\chi_j\subset\Gamma_j$.

Тоді, за означенням, простір $H^{s,\varphi}(\Gamma)$ складається з усіх розподілів $h$ на $\Gamma$ таких, що
$$
(\chi_{j}h)\circ\pi_{j}\in H^{s,\varphi}(\mathbb{R}^{n-1})
$$
для кожного номера $j\in\{1,\ldots,\lambda\}$. Тут, звісно,
$(\chi_{j}h)\circ\pi_{j}$ є зображенням розподілу $h$ у локальній карті $\pi_{j}$. У~просторі $H^{s,\varphi}(\Gamma)$ уведено норму за формулою
$$
\|h\|_{H^{s,\varphi}(\Gamma)}:=\biggl(\,\sum_{j=1}^{\lambda}
\|(\chi_{j}h)\circ\pi_{j}\|_{H^{s,\varphi}(\mathbb{R}^{n-1})}^{2}
\biggr)^{1/2}.
$$
Цей простір гільбертів і сепарабельний відносно цієї норми. Важливо, що він з точністю до еквівалентності норм не залежить від зазначеного вибору атласу і розбиття одиниці \cite[теорема~2.3]{MikhailetsMurach14}. Множина $C^{\infty}(\Gamma)$ щільна у просторі $H^{s,\varphi}(\Gamma)$.

Уведені гільбертові функціональні простори утворюють уточнені соболєвські шкали
\begin{equation}\label{1f8}
\{H^{s,\varphi}(\Omega):\,s\in\mathbb{R},\varphi\in\mathcal{M}\}
\quad\mbox{і}\quad
\{H^{s,\varphi}(\Gamma):\,s\in\mathbb{R},\varphi\in\mathcal{M}\}
\end{equation}
на $\Omega$ і $\Gamma$ відповідно. Вони містять гільбертові соболєвські шкали: якщо $\varphi(\cdot)\equiv1$, то $H^{s,\varphi}(\Omega)=:H^{s}(\Omega)$ і $H^{s,\varphi}(\Gamma)=:H^{s}(\Gamma)$ є простори Соболєва порядку
$s\in\mathbb{R}$.

Для шкал \eqref{1f8} виконуються компактні і щільні вкладення \eqref{1f7}, якщо у формулі \eqref{1f7} замінити $\mathbb{R}^{n}$ на $\Omega$ або $\Gamma$ відповідно.

\section{Основні результати}\label{1sec4}

У цьому розділі сформулюємо основні результати статті про властивості еліптичної крайової задачі \eqref{1f1}, \eqref{1f2} в уточненій соболєвській шкалі.

Пов'яжемо із задачею \eqref{1f1}, \eqref{1f2} гільбертові простори
\begin{equation*}
\mathcal{D}^{s,\varphi}(\Omega,\Gamma):=H^{s,\varphi}(\Omega)\oplus
\bigoplus_{k=1}^{\varkappa}H^{s+r_{k}-1/2,\varphi}(\Gamma),
\end{equation*}
та
\begin{equation*}
\mathcal{E}^{s-2q,\varphi}(\Omega,\Gamma):=H^{s-2q,\varphi}(\Omega)\oplus
\bigoplus_{j=1}^{q+\varkappa}H^{s-m_{j}-1/2,\varphi}(\Gamma),
\end{equation*}
де $s\in\mathbb{R}$ і $\varphi\in\mathcal{M}$. У соболєвському випадку, коли $\varphi(\cdot)\equiv1$, будемо пропускати індекс $\varphi$ в позначеннях просторів $\mathcal{D}^{s,\varphi}(\Omega,\Gamma)$ і $\mathcal{E}^{s-2q,\varphi}(\Omega,\Gamma)$.

Позначимо через $\mathcal{N}$ лінійний простір усіх розв'язків
$$
(u,v_{1},\ldots,v_{\varkappa})\in C^{\infty}(\overline{\Omega})\times (C^{\infty}(\Gamma))^{\varkappa}
$$
крайової задачі \eqref{1f1}, \eqref{1f2} в однорідному випадку, коли $f=0$ в $\Omega$ і кожне $g_{j}=0$ на~$\Gamma$. Аналогічно, позначимо через $\mathcal{N}_{\star}$ лінійний простір усіх розв'язків
$$
(\omega,w_{1},\ldots,w_{m-2q+1},h_{1},\ldots,h_{q+\varkappa})\in C^{\infty}(\overline{\Omega})\times(C^{\infty}(\Gamma))^{m-q+1+\varkappa}
$$
формально спряженої крайової задачі \eqref{1f4}\,--\,\eqref{1f6} в однорідному випадку, коли $\eta=0$ в $\Omega$ і кожне $\psi_{k}=0$ на~$\Gamma$. Оскільки обидві задачі еліптичні за Лавруком в $\Omega$, то простори $\mathcal{N}$ і $\mathcal{N}_{\star}$ скінченновимірні \cite[наслідок~4.1.1]{KozlovMazyaRossmann97}.

\begin{theorem}\label{1th1}
Нехай $s>m+1/2$ і $\varphi\in\mathcal{M}$. Тоді відображення \eqref{1f3} продовжується єдиним чином (за неперервністю) до обмеженого оператора
\begin{equation}\label{1f9}
\Lambda:\mathcal{D}^{s,\varphi}(\Omega,\Gamma)\rightarrow
\mathcal{E}^{s-2q,\varphi}(\Omega,\Gamma).
\end{equation}
Цей оператор нетерів. Його ядро збігається з простором $\mathcal{N}$, а область значення складається з усіх векторів
\begin{equation}\label{1f10}
(f,g):=
(f,g_{1},\ldots,g_{q+\varkappa})\in
\mathcal{E}^{s-2q,\varphi}(\Omega,\Gamma)
\end{equation}
таких, що
\begin{equation}\label{1f11}
\begin{gathered}
(f,\omega)_\Omega+\sum_{j=1}^{m-2q+1}(D_{\nu}^{j-1}f,w_{j})_{\Gamma}+
\sum_{j=1}^{q+\varkappa}(g_j,h_{j})_{\Gamma}=0 \\
\mbox{для всіх} \quad (\omega,w_{1},\ldots,w_{m-2q+1},h_{1},\ldots,h_{q+\varkappa})\in \mathcal{N}_{\star}.
\end{gathered}
\end{equation}
Індекс оператора \eqref {1f9} дорівнює $\dim\mathcal{N}-\dim\mathcal{N}_{\star}$ й не залежить від $s$ та~$\varphi$.
\end{theorem}

Зауважимо, що в умові \eqref{1f11} функції
$$
D_{\nu}^{j-1}f\in H^{s-2q-j+1/2}(\Gamma)\subset L_{2}(\Gamma)
$$
коректно означені для кожного номера $j\in\{1,\ldots,m-2q+1\}$ на підставі \cite[теорема 4.13(ii)]{MikhailetsMurach14} і нерівності $s>m+1/2$.

Стосовно цієї теореми нагадаємо, що лінійний обмежений оператор $T:X\rightarrow Y$, де $X,Y$~--- банахові простори, називають нетеровим, якщо його ядро $\ker T$ і коядро $Y/T(X)$ скінченновимірні. Нетерів оператор має замкнену область значень $T(X)$ (див., наприклад, \cite[лемма~19.1.1]{Hermander87}) і скінченний індекс
$$
\mathrm{ind}\,T:=\dim\ker T-\dim(Y/\,T(X)).
$$

Зауважимо, що умову $s>m+1/2$ в теоремі~\ref{1th1} не можна відкинути або послабити. Зокрема, якщо $s\leq m_{j}+1/2$ і $\varphi(\cdot)\equiv1$, то відображення $u\mapsto B_{j}u$, де $u\in C^{\infty}(\overline{\Omega})$, не можна продовжити до неперервного лінійного оператора, що діє з простору Соболєва $H^{s}(\Omega)$ у лінійний топологічний простір $\mathcal{D}'(\Gamma)$ усіх розподілів на~$\Gamma$ (див., наприклад, \cite[зауваження~3.5]{MikhailetsMurach14}).

Якщо $\mathcal{N}=\{0\}$ і $\mathcal{N}_{\star}=\{0\}$, то оператор \eqref{1f9} є ізоморфізмом простору $\mathcal{D}^{s,\varphi}(\Omega,\Gamma)$ на простір $\mathcal{E}^{s-2q,\varphi}(\Omega,\Gamma)$. У загальній ситуації цей оператор породжує ізоморфізм між деякими їх (замкненими) підпросторами скінченної ковимірності. Ці підпростори і проектори на них будуємо у такий спосіб.

Розглянемо розклад простору $\mathcal{D}^{s,\varphi}(\Omega,\Gamma)$, де $s>m+1/2$ і $\varphi\in\mathcal{M}$, у таку пряму суму його  підпросторів:
\begin{equation}\label{1f12}
\begin{gathered}
\mathcal{D}^{s,\varphi}(\Omega,\Gamma)=\mathcal{N}\dotplus
\biggl\{(u,v_{1},\ldots,v_{\varkappa})\in
\mathcal{D}^{s,\varphi}(\Omega,\Gamma):\\
(u,u^{\circ})_\Omega+
\sum_{k=1}^{\varkappa}(v_{k},v_{k}^{\circ})_{\Gamma}\;\,
\mbox{для кожного}\;\,
(u^{\circ},v_{1}^{\circ},\ldots,v_{\varkappa}^{\circ})
\in\mathcal{N}\biggr\}.
\end{gathered}
\end{equation}
Таке зображення існує, оскільки воно є звуженням розкладу простору $L_{2}(\Omega)\oplus(L_{2}(\Gamma))^{\varkappa}$ в ортогональну суму підпростору $\mathcal{N}$ та його доповнення. Тут
\begin{equation}\label{embedding-in-L_2}
\mathcal{D}^{s,\varphi}(\Omega,\Gamma)\hookrightarrow L_{2}(\Omega)\oplus(L_{2}(\Gamma))^{\varkappa}
\quad\mbox{при}\quad s>m+1/2
\end{equation}
на підставі умови \eqref{nat-assump}.

Стосовно розкладу простору $\mathcal{E}^{s-2q,\varphi}(\Omega,\Gamma)$ скористаємося таким результатом.

\begin{lemma}\label{1lem1}
Існує скінченновимірний простір
$$
\mathcal{G}\subset C^{\infty}(\overline{\Omega})\times(C^{\infty}(\Gamma))^{q+\varkappa}
$$
такий, що для довільних $s>m+1/2$ і $\varphi\in\mathcal {M}$ є правильним наступний розклад простору $\mathcal{E}^{s-2q,\varphi}(\Omega,\Gamma)$ у пряму суму його підпросторів:
\begin{equation}\label{1f13}
\begin{gathered}
\mathcal{E}^{s-2q,\varphi}(\Omega,\Gamma)=\mathcal{G}\dotplus\\
\bigl\{(f,g)\in
\mathcal{E}^{s-2q,\varphi}(\Omega,\Gamma):
\mbox{виконується \eqref{1f11}}\bigr\}.
\end{gathered}
\end{equation}
При цьому $\dim\mathcal{G}=\dim\mathcal{N}_{\star}$.
\end{lemma}

Позначимо через $\mathcal{P}$ і $\mathcal{Q}$ відповідно проектори просторів $\mathcal{D}^{s,\varphi}(\Omega,\Gamma)$ і $\mathcal{E}^{s-2q,\varphi}(\Omega,\Gamma)$ на другий доданок у сумах \eqref{1f12} і \eqref{1f13} паралельно першому доданку. Ці проектори не залежать (як відо\-браження) від $s$ і $\varphi$.

\begin{theorem}\label{1th2}
Нехай $s>m+1/2$ і $\varphi\in\mathcal{M}$. Тоді звуження відображення \eqref{1f9} на підпростір $\mathcal{P}(\mathcal{D}^{s,\varphi}(\Omega,\Gamma))$ є ізоморфізмом
\begin{equation}\label{1f14}
\Lambda:\,\mathcal{P}(\mathcal\mathcal{D}^{s,\varphi}(\Omega,\Gamma))
\leftrightarrow
\mathcal{Q}(\mathcal{E}^{s-2q,\varphi}(\Omega,\Gamma)).
\end{equation}
\end{theorem}

Дослідимо властивості узагальнених розв'язків еліптичної крайової задачі \eqref{1f1}, \eqref{1f2}. Попередньо дамо означення такого розв'язку.
Покладемо
$$
\mathcal{D}^{m+1/2+}(\Omega,\Gamma):=
\bigcup_{\substack{s>m+1/2,\\\varphi\in\mathcal{M}}}
\mathcal{D}^{s,\varphi}(\Omega,\Gamma)=\bigcup_{s>m+1/2}\mathcal{D}^{s}(\Omega,\Gamma);
$$
тут остання рівність правильна з огляду на властивість \eqref{1f7}.
Вектор
\begin{equation}\label{generalized-solution}
(u,v):=(u,v_{1},\ldots,v_{\varkappa})\in\mathcal{D}^{m+1/2+}(\Omega,\Gamma)
\end{equation}
називаємо (сильним) узагальненим розв'язком крайової задачі \eqref{1f1}, \eqref{1f2} з правої частиною
$$
(f,g):=(f,g_{1},\ldots,g_{q+\varkappa})\in \mathcal{D}'(\Omega)\times\bigl(\mathcal{D}'(\Gamma)\bigr)^{q+\varkappa},
$$
якщо $\Lambda(u,v)=(f,g)$, де $\Lambda$~--- оператор \eqref{1f9} для деяких параметрів $s>m+1/2$ і $\varphi\in\mathcal{M}$. Тут $\mathcal{D}'(\Omega)$~--- лінійний топологічний простір усіх розподілів, заданих в області $\Omega$. Наведене означення узагальненого розв'язку коректне, оскільки не залежить від $s$ і~$\varphi$.

\begin{theorem}\label{1th3}
Нехай $\varphi\in\mathcal{M}$, а числа $s,\lambda\in\mathbb{R}$   задовольняють нерівності $s>m+1/2$ і $0<\lambda<s-m+1/2$. Нехай також функції $\chi,\eta\in C^{\infty}(\overline{\Omega})$ задовольняють умову $\eta =1$ в околі $\mathrm{supp}\,\chi$. Тоді існує число $c=c(s,\varphi,\lambda,\chi,\eta)>0$ таке, що
\begin{equation}\label{1f15}
\|\chi(u,v)\|_{\mathcal{D}^{s,\varphi}(\Omega,\Gamma)}\leq c\,\bigl(\,\|\eta\Lambda(u,v)\|_{\mathcal{E}^{s-2q,\varphi}(\Omega,\Gamma)}+
\|\eta(u,v)\|_{\mathcal{D}^{s-\lambda,\varphi}(\Omega,\Gamma)}\bigl)
\end{equation}
для довільного вектора $(u,v)\in\mathcal{D}^{s,\varphi}(\Omega,\Gamma)$. Тут $c$ не залежить від~$(u,v)$.
\end{theorem}

Тут, звісно,
$$
\chi(u,v):=(\chi u,(\chi\!\upharpoonright\!\Gamma) v_{1},\ldots,(\chi\!\upharpoonright\!\Gamma)v_{\varkappa})
$$
і аналогічно розуміємо вираз $\eta\Lambda(u,v)$.

\begin{remark}\label{rem1}\rm
Якщо $\chi=\eta=1$, то нерівність \eqref{1f15} є глобальною апріорною оцінкою узагальненого розв'язку $u$ еліптичної крайової задачі \eqref{1f1}, \eqref{1f2}. У цьому випадку умову $\lambda<s-m+1/2$ можна прибрати.  Взагалі, нерівність \eqref{1f15} є локальною апріорною оцінкою розв'язку~$u$. Справді, для кожної непорожньої відкритої (у топології $\overline{\Omega}$) підмножини множини $\overline{\Omega}$, можна вибрати функції $\chi,\eta$ так, щоб вони задовольняли умову теореми~\ref{1th3} і їх носії лежали в цій підмножині. Якщо $0<\lambda\leq1$, то у нерівності \eqref{1f15} можна узяти $\chi\Lambda(u,v)$ замість $\eta\Lambda(u,v)$.
\end{remark}

Дослідимо регулярність узагальнених розв'язків еліптичної крайової задачі \eqref{1f1}, \eqref{1f2}. Нехай $V$~--- довільна відкрита підмножина простору $\mathbb{R}^{n}$ така, що
$\Omega_0:=\Omega\cap V\neq\varnothing$. Покладемо $\Gamma_{0}:=\Gamma\cap V$ (можливий випадок, коли $\Gamma_{0}=\varnothing$). Для довільних параметрів $\sigma\in\mathbb{R}$ і $\varphi\in\mathcal{M}$ введемо локальні аналоги просторів Хермандера $H^{\sigma,\varphi}(\Omega)$ і $H^{\sigma,\varphi}(\Gamma)$. А саме,
\begin{gather*}
H^{\sigma,\varphi}_{\mathrm{loc}}(\Omega_{0},\Gamma_{0}):=\bigl\{u\in \mathcal{D}'(\Omega):\,\chi u\in H^{\sigma,\varphi}(\Omega)\\
\mbox{для всіх}\;\,\chi\in C^{\infty}(\overline{\Omega})\;\,\mbox{таких, що}\;\,\mathrm{supp}\,\chi\subset\Omega_0\cup\Gamma_{0}\bigr\}.
\end{gather*}
Аналогічно покладемо
\begin{gather*}
H^{\sigma,\varphi}_{\mathrm{loc}}(\Gamma_{0}):=\bigl\{h\in \mathcal{D}'(\Gamma):\,\chi h\in H^{\sigma,\varphi}(\Gamma)\\
\mbox{для всіх}\;\,\chi\in C^{\infty}(\Gamma)\;\,\mbox{таких, що}\;\,
\mathrm{supp}\,\chi\subset\Gamma_{0}\bigr\}.
\end{gather*}
Позначимо
\begin{gather*}
\mathcal{D}^{s,\varphi}_{\mathrm{loc}}(\Omega_{0},\Gamma_{0}):=
H^{s,\varphi}_{\mathrm{loc}}(\Omega_{0},\Gamma_{0})
\times\prod\limits_{k=1}^{\varkappa}
H^{s+r_k-1/2,\varphi}_{\mathrm{loc}}(\Gamma_{0}),\\
\mathcal{E}^{s-2q,\varphi}_{\mathrm{loc}}(\Omega_{0},\Gamma_{0}):=
H^{s-2q,\varphi}_{\mathrm{loc}}(\Omega_{0},\Gamma_{0})
\times\prod\limits_{j=1}^{q+\varkappa}
H^{s-m_j-1/2,\varphi}_{\mathrm{loc}}(\Gamma_{0}),
\end{gather*}
де $s\in\mathbb{R}$ і $\varphi\in\mathcal{M}$.

\begin{theorem}\label{1th4}
Нехай вектор $(u,v)\in\mathcal{D}^{m+1/2+}(\Omega,\Gamma)$ є узагальненим розв'язком еліптичної крайової задачі \eqref{1f1}, \eqref{1f2}, праві частини якої задовольняють умову $(f,g)\in\mathcal{E}^{s-2q,\varphi}_{\mathrm{loc}}(\Omega_{0},\Gamma_{0})$
для деяких параметрів $s>m+1/2$ і $\varphi\in\mathcal{M}$. Тоді $(u,v)\in\mathcal{D}^{s,\varphi}_{\mathrm{loc}}(\Omega_{0},\Gamma_{0})$.
\end{theorem}

Відмітимо важливі окремі випадки цієї теореми. Якщо $\Omega_{0}=\Omega$ і $\Gamma_{0}=\Gamma$, то локальні простори $\mathcal{D}^{s,\varphi}_{\mathrm{loc}}(\Omega_{0},\Gamma_{0})$ i $\mathcal{E}^{s-2q,\varphi}_{\mathrm{loc}}(\Omega_{0},\Gamma_{0})$ збігаються з просторами $\mathcal{D}^{s,\varphi}(\Omega,\Gamma)$ і $\mathcal{E}^{s-2q,\varphi}(\Omega,\Gamma)$ відповідно. Тому теорема~\ref{1th4} стверджує, що регулярність узагальненого розв'язку $(u,v)$ підвищується глобально, тобто в усій області $\Omega$ аж до її межі $\Gamma$. Якщо  $\Gamma_{0}=\varnothing$ і $\Omega_{0}=\Omega$, то згідно з цією теоремою регулярність розв'язку $(u,v)$ підвищується в околах усіх внутрішніх точок замкненої області~$\overline\Omega$.

У соболєвському випадку, коли $\varphi(\cdot)\equiv1$, теореми \ref{1th1}\,--\,\ref{1th4} або їх версії відомі, див., наприклад,  монографію В.~О.~Козлова, В.~Г.~Маз'ї і Й.~Россмана \cite[п.~4.1]{KozlovMazyaRossmann97}, статті І.~Я.~Ройтберг \cite{RoitbergInna97, RoitbergInna98} та монографію А.~Я.~Ройтберга \cite[п.~2.4]{Roitberg99}.

\section{Застосування}\label{1sec5}

Розглянемо застосування уточненої соболєвської шкали до питання про
достатні умови неперервності узагальнених похідних (заданого порядку) компонент розв'язку \eqref{generalized-solution} еліптичної крайової задачі \eqref{1f1}, \eqref{1f2}. Ці умови випливають з теореми~\ref{1th4} і такого наслідку з теореми вкладення Хермандера \cite[теорема~2.2.7]{Hermander65}:

\begin{proposition}\label{1prop1}
Нехай $0\leq l\in\mathbb{Z}$ і $\varphi\in\mathcal{M}$. Тоді кожне з вкладень  $H^{l+n/2,\varphi}(\Omega)\subset C^{l}(\overline{\Omega})$ і
$H^{l+(n-1)/2,\varphi}(\Gamma)\subset C^{l}(\Gamma)$ еквівалентне виконанню умови
\begin{equation}\label{1f45}
\int\limits_{1}^{\infty}\frac{dt}{t\,\varphi^{2}(t)}<\infty.
\end{equation}
Ці вкладення компактні.
\end{proposition}

Твердження \ref{1prop1} обґрунтовано в \cite[теореми~2.8 і 3.4]{MikhailetsMurach14}. Можна сказати, що воно є уточненням відомої теореми вкладення Соболєва, згідно з якою
\begin{gather*}
s>l+n/2\;\Leftrightarrow\;
H^{s}(\Omega)\subset C^{l}(\overline{\Omega}),\\
s>l+(n-1)/2\;\Leftrightarrow\;H^{s}(\Gamma)\subset C^{l}(\Gamma),
\end{gather*}
де $0\leq l\in\mathbb{Z}$.

Сформулюємо достатні умови, згадані на початку цього пункту.

\begin{theorem}\label{1th5} Нехай ціле число $l\geq0$. Припустимо, що вектор $(u,v)\in\mathcal{D}^{m+1/2+}(\Omega,\Gamma)$ є узагальненим розв'язком еліптичної крайової задачі \eqref{1f1}, \eqref{1f2}, де
\begin{equation}\label{1th5-cond}
(f,g)\in
\mathcal{E}^{l+n/2-2q,\varphi}_{\mathrm{loc}}(\Omega_{0},\Gamma_{0})
\end{equation}
для деякого функціонального параметра $\varphi\in\mathcal{M}$, який задовольняє умову \eqref{1f45}. Тоді $u\in C^{l}(\Omega_{0}\cup\Gamma_{0})$.
\end{theorem}

\begin{theorem}\label{1th6}
Нехай задано цілі числа $l\geq0$ і $k\in\{1,...,\varkappa\}$. Припустимо, що $\Gamma_{0}\neq\varnothing$ і вектор $(u,v)\in\mathcal{D}^{m+1/2+}(\Omega,\Gamma)$ є узагальненим розв'язком еліптичної крайової задачі \eqref{1f1}, \eqref{1f2}, де
\begin{equation}\label{1th6-cond}
(f,g)\in
\mathcal{E}^{l-r_k+n/2-2q,\varphi}_{\mathrm{loc}}(\Omega_{0},\Gamma_{0})
\end{equation}
для деякого функціонального параметра $\varphi\in\mathcal{M}$, який задовольняє умову \eqref{1f45}. Тоді $v_{k}\in C^{l}(\Gamma_{0})$.
\end{theorem}

\begin{remark}\label{rem2} \rm
Умова \eqref{1f45} є точною у теоремах \ref{1th5} і \ref{1th6}, якщо число $l$ достатньо велике. А саме, нехай $0\leq l\in\mathbb{Z}$ і $\varphi\in\mathcal{M}$. Припустимо, що $l+n/2>m+1/2$; тоді умова \eqref{1f45} еквівалентна імплікації
\begin{equation}\label{rem1-impl-a}
\bigl((u,v)\in\mathcal{D}^{m+1/2+}(\Omega,\Gamma)\;\mbox{задовольняє умову теореми~$\ref{1th5}$}\bigr)\Rightarrow
u\in C^{l}(\Omega_{0}\cup\Gamma_{0}).
\end{equation}
Окрім того, нехай $k\in\{1,...,\varkappa\}$ і $\Gamma_{0}\neq\varnothing$. Припустимо, що $l-r_{k}+n/2>m+1/2$; тоді умова \eqref{1f45} еквівалентна імплікації
\begin{equation}\label{rem1-impl-b}
\bigl((u,v)\in\mathcal{D}^{m+1/2+}(\Omega,\Gamma)\;\mbox{задовольняє умову теореми~$\ref{1th6}$}\bigr)\Rightarrow 
v_{k}\in C^{l}(\Gamma_{0}).
\end{equation}
\end{remark}

Якщо б ми сформулювати версії теорем \ref{1th5} і \ref{1th6} для просторів Соболєва, то замість умов \eqref{1th5-cond} і \eqref{1th6-cond} використали б таку більш сильну умову: $(f,g)\in\mathcal{H}^{s-2q,1}_{\mathrm{loc}}(\Omega_{0},\Gamma_{0})$
для деякого $s>l+n/2-2q$ у випадку теореми~\ref{1th5} або для деякого
$s>l-r_k+n/2-2q$ у випадку теореми~\ref{1th6}.

Сформулюємо також достатню умову того, що узагальнений розв'язок $(u,v)$ крайової задачі \eqref{1f1}, \eqref{1f2} є класичним, тобто $u\in\nobreak C^{2q}(\Omega)\cap C^{m}(U_{\sigma}\cup \Gamma)$ для деякого числа $\sigma>0$ і $v_{k}\in C^{m+r_{k}}(\Gamma)$ для кожного номера $k\in\{1,...,\varkappa\}$. Тут
$$
U_{\sigma}:=\{x\in\Omega:\mathrm{dist}(x,\Gamma)<\sigma\}.
$$
Якщо розв'язок $(u,v)$ цієї задачі є класичним, то її ліві частини обчислюються за допомогою класичних похідних і є неперервними функціями на $\Omega$ і $\Gamma$ відповідно.

\begin{theorem}\label{1th7} Нехай вектор $(u,v)\in\mathcal{D}^{m+1/2+}(\Omega,\Gamma)$ є узагальненим розв'язком еліптичної крайової задачі \eqref{1f1}, \eqref{1f2}, де
\begin{gather}\label{1f46}
f\in H^{n/2,\varphi_1}_{\mathrm{loc}}(\Omega,\varnothing)\cap H^{m+n/2-2q,\varphi_2}_{\mathrm{loc}}(U_{\sigma},\Gamma),\\
g_j\in H^{m-m_j+(n-1)/2,\varphi_2}(\Gamma)
\quad\mbox{при}\quad j\in\{1,\ldots,q+\varkappa\} \label{1f47}
\end{gather}
для деякого числа $\sigma>0$ і деяких функцій $\varphi_1,\varphi_2\in\mathcal{M}$, які задовольняють умову \eqref{1f45} при $\varphi:=\varphi_{1}$ і $\varphi:=\varphi_{2}$ відповідно. Тоді розв'язок $(u,v)$ класичний.
\end{theorem}

\section{Доведення}\label{1sec6}

Доведемо теореми \ref{1th1}\,--\,\ref{1th7} і лему~\ref{1lem1} та обґрунтуємо зауваження \ref{rem1} і~\ref{rem2}.

\medskip

\nobreak\textbf{Доведення теореми~\ref{1th1}.}  У  випадку просторів Соболєва, коли $\varphi(\cdot)\equiv1$ і $s>m+1/2$, ця теорема (для загальних еліптичних систем) доведена в статтях \cite{RoitbergInna97, RoitbergInna98} та в монографії \cite[теорема 2.4.1]{Roitberg99} за виключенням вказаного зв'язку скінченновимірного простору $\mathcal{N}_{\star}$ з формально спряженою задачею \eqref{1f4}\,--\,\eqref{1f6}.
За додаткового припущення $s\in\mathbb{Z}$, у повному обсязі теорема~\ref{1th1} міститься у результаті, встановленому в монографії В.~О.~Козлова, В.~Г.~Маз'ї і Й.~Россмана \cite[наслідок 4.1.1]{KozlovMazyaRossmann97}. Доведемо, що і для кожного  дробового $s>m+1/2$ висновок цієї теореми правильний у повному обсязі.

Згідно з \cite[теорема~2.4.1]{Roitberg99} відображення \eqref{1f3} продовжується за неперервністю до обмеженого і нетерового оператора
\begin{equation}\label{roit--oper}
\begin{gathered}
\Lambda:H^{s,(m+1)}(\Omega)\oplus
\bigoplus_{k=1}^{\varkappa}H^{s+r_{k}-1/2}(\Gamma)\rightarrow\\
\rightarrow H^{s-2q,(m+1-2q)}(\Omega)\oplus
\bigoplus_{j=1}^{q+\varkappa}H^{s-m_{j}-1/2}(\Gamma)
\end{gathered}
\end{equation}
для кожного $s\in\mathbb{R}$. Тут $H^{s,(r)}(\Omega)$, де $s\in\mathbb{R}$ і $1\leq r\in\mathbb{Z}$, є гільбертовим простором Соболєва\,--\,Ройтберга \cite[п.~2.1]{Roitberg96}. Якщо $s\notin\{1/2,\ldots,r-1/2\}$, то, за означенням, $H^{s,(r)}(\Omega)$~--- поповнення простору $C^{\infty}(\overline{\Omega})$ за нормою
$$
\|u\|_{H^{s,(r)}(\Omega)}:=
\biggl(\|u\|_{H^{s,(0)}(\Omega)}^{2}+
\sum_{k=1}^{r}\;\|(D_{\nu}^{k-1}u)\!\upharpoonright\!\Gamma\|
_{H^{s-k+1/2}(\Gamma)}^{2}\biggr)^{1/2}.
$$
Тут $H^{s,(0)}(\Omega):=H^{s}(\Omega)$ для $s\geq0$; якщо $s<0$, то $H^{s,(0)}(\Omega)$ є дуальний гільбертів простір до простору $H^{-s}(\Omega)$ відносно розширення за неперервністю скалярного добутку в $L_{2}(\Omega)$.

Зауважимо, що для просторів Соболєва\,--\,Ройтберга виконується неперервне і щільне вкладення $H^{s+\delta,(r)}(\Omega)\hookrightarrow H^{s,(r)}(\Omega)$ при $\delta>0$. Окрім того, якщо $s>r-1/2$, то простори $H^{s,(r)}(\Omega)$ і $H^{s}(\Omega)$ рівні як поповнення лінійного многовиду $C^{\infty}(\overline{\Omega})$ за еквівалентними нормами. Тому оператор \eqref{1f9}, де $\varphi(\cdot)\equiv1$, і оператор \eqref{roit--oper} рівні при $s>m+1/2$.

З огляду на формулу \eqref{roit--oper} покладемо
\begin{gather*}
\mathcal{D}^{s,(m+1)}(\Omega,\Gamma):=H^{s,(m+1)}(\Omega)\oplus
\bigoplus_{k=1}^{\varkappa}H^{s+r_{k}-1/2}(\Gamma),\\
\mathcal{E}^{s-2q,(m+1-2q)}(\Omega,\Gamma):=
H^{s-2q,(m+1-2q)}(\Omega)\oplus
\bigoplus_{j=1}^{q+\varkappa}H^{s-m_{j}-1/2}(\Gamma).
\end{gather*}

Згідно із згаданим результатом \cite[теорема~2.4.1]{Roitberg99} ядро оператора \eqref{roit--oper} дорівнює $\mathcal{N}$, а область значень складається з усіх векторів
$$
(f,g)\in\mathcal{E}^{s-2q,(m+1-2q)}(\Omega,\Gamma),
$$
які задовольняють  умову \eqref{1f11}, у котрій замість $\mathcal{N}_{\star}$ фігурує деякий скінченновимірний простір, що лежить в
$$
C^{\infty}(\overline{\Omega})\times(C^{\infty}(\Gamma))^{m-q+1+\varkappa}
$$
і не залежить від $s$. Звідси випливає рівність
\begin{equation*}
\Lambda(\mathcal{D}^{s_{2},(m+1)}(\Omega,\Gamma))=
\mathcal{E}^{s_{2}-2q,(m+1-2q)}(\Omega,\Gamma)\cap
\Lambda(\mathcal{D}^{s_{1},(m+1)}(\Omega,\Gamma))
\end{equation*}
для довільних чисел $s_{1},s_{2}\in\mathbb{R}$ таких, що $s_{1}<s_{2}$.
Зокрема,
\begin{gather*}
\Lambda(\mathcal{D}^{s}(\Omega,\Gamma))=
\mathcal{E}^{s-2q}(\Omega,\Gamma)\cap
\Lambda(\mathcal{D}^{m,(m+1)}(\Omega,\Gamma)),\\
\quad\mbox{якщо}\quad m+1/2<s\in\mathbb{R}.
\end{gather*}
На підставі \cite[теорема 4.1.4]{KozlovMazyaRossmann97} простір $\Lambda(\mathcal{D}^{m,(m+1)}(\Omega,\Gamma))$ складається з усіх векторів $(f,g)\in\mathcal{E}^{m-2q,(m+1-2q)}(\Omega,\Gamma)$, які задовольняють умову \eqref{1f11}. Тому для довільного дійсного $s>m+1/2$ область значень $\Lambda(\mathcal{D}^{s}(\Omega,\Gamma))$ оператора \eqref{1f9}, де $\varphi(\cdot)\equiv1$, є такою як це стверджується у теоремі~\ref{1th1}. Таким чином, у соболєвському випадку ця теорема доведена.

Теорему~\ref{1th1} для довільного $\varphi\in\mathcal{M}$ виведемо тепер із соболєвського випадку за допомогою інтерполяції з функціональним параметром (її означення і властивості наведено, наприклад, у монографії \cite[п. 1.1]{MikhailetsMurach14}). Нехай $s>m+1/2$ і $\varphi\in\mathcal{M}$. Покладемо $\varepsilon:=(s-m-1/2)/2>0$. Відображення \eqref{1f3} продовжується за неперервністю до обмежених і нетерових операторів
\begin{equation}\label{1f17}
\begin{gathered}
\Lambda:\mathcal{D}^{s\mp\varepsilon}(\Omega,\Gamma)\rightarrow
\mathcal{E}^{s\mp\varepsilon-2q}(\Omega,\Gamma),
\end{gathered}
\end{equation}
що діють у парах соболєвських просторів. Вони мають спільне ядро $\mathcal{N}$ та однаковий індекс, рівний $\dim \mathcal{N}-\dim \mathcal{N}_{\star}$. Окрім того,
\begin{equation}\label{1f18}
\Lambda(\mathcal{D}^{s\mp\varepsilon}(\Omega,\Gamma))
=\bigl\{(f,g)\in\mathcal{E}^{s\mp\varepsilon-2q}(\Omega,\Gamma):
\mbox{виконується}\;\eqref{1f11}\bigr\}.
\end{equation}
Означимо інтерполяційний параметр $\psi$ за формулами $\psi(t):=t^{1/2}
\varphi(t^{1/(2\varepsilon)})$ при $t\geq1$ і $\psi(t):=\varphi(1)$ при $0<t<1$. Застосувавши до \eqref{1f17} інтерполяцію з параметром $\psi$ та використавши теорему про інтерполяцію нетерових операторів \cite[теорема 1.7]{MikhailetsMurach14}, отримаємо нетерів обмежений оператор
\begin{equation}\label{1f19}
\Lambda:[\mathcal{D}^{s-\varepsilon}(\Omega,\Gamma),
\mathcal{D}^{s+\varepsilon}(\Omega,\Gamma)]_{\psi}\rightarrow
[\mathcal{E}^{s-\varepsilon-2q}(\Omega,\Gamma),
\mathcal{E}^{s+\varepsilon-2q,\varphi}(\Omega,\Gamma)]_{\psi}.
\end{equation}
(Тут $[H_{0},H_{1}]_{\psi}$ позначає гільбертів простір, отриманий інтерполяцією з параметром $\psi$ допустимої пари гільбертових просторів $H_{0}$ і $H_{1}$.) Оператор \eqref{1f19} є звуженням відображення \eqref{1f17}, заданого на $\mathcal{D}^{s-\varepsilon}(\Omega,\Gamma)$. Отже, він є продовженням за неперервністю відображення
\eqref{1f3}.

Опишемо інтерполяційні простори, в яких діє оператор \eqref{1f19}. На підставі інтерполяційних теорем 1.5, 2.2 і 3.2 з монографії \cite{MikhailetsMurach14} маємо такі рівності просторів разом з еквівалентністю норм у них:
\begin{gather*}
[\mathcal{D}^{s-\varepsilon}(\Omega,\Gamma),
\mathcal{D}^{s+\varepsilon}(\Omega,\Gamma)]_{\psi}=
[H^{s-\varepsilon}(\Omega),H^{s+\varepsilon}(\Omega)]_{\psi}\oplus\\
\oplus\bigoplus_{k=1}^{\varkappa}\,
[H^{s+r_{k}-1/2-\varepsilon}(\Gamma),
H^{s+r_{k}-1/2+\varepsilon}(\Gamma)]_{\psi}=\\
=H^{s,\varphi}(\Omega)\oplus
\bigoplus_{k=1}^{\varkappa}H^{s+r_{k}-1/2,\varphi}(\Gamma)=
\mathcal{D}^{s,\varphi}(\Omega,\Gamma).
\end{gather*}
Аналогічно
$$
[\mathcal{E}^{s-2q-\varepsilon}(\Omega,\Gamma),
\mathcal{E}^{s-2q+\varepsilon}(\Omega,\Gamma)]_{\psi}=
\mathcal{E}^{s-2q,\varphi}(\Omega,\Gamma).
$$
Таким чином, обмежений нетерів оператор \eqref{1f19} є оператором \eqref{1f9} з теореми~\ref{1th1}.

Згідно з \cite[теорема 1.7]{MikhailetsMurach14} ядро оператора \eqref{1f9} і його індекс
збігаються відповідно зі спільним ядром $\mathcal{N}$ і однаковим індексом $\dim \mathcal{N}-\dim \mathcal{N}_{\star}$
операторів \eqref{1f17}. Окрім того, на підставі рівності \eqref{1f18} робимо висновок, що область значень оператора \eqref{1f9} дорівнює
\begin{equation*}
\mathcal{E}^{s-2q,\varphi}(\Omega,\Gamma)\cap
\Lambda(\mathcal{D}^{s-\varepsilon}(\Omega,\Gamma))
=\bigl\{(f,g)\in\mathcal{E}^{s-2q,\varphi}(\Omega,\Gamma):\,
\mbox{виконується}\;\eqref{1f11}\bigr\}.
\end{equation*}
Отже, доведено всі властивості оператора \eqref{1f9}, зазначені в теоремі~\ref{1th1}.

Теорему~\ref{1th1} доведено.

\medskip

\noindent\textbf{Доведення леми~\ref{1lem1}.} Скористаємося нетеровим обмеженим оператором \eqref{roit--oper} у випадку, коли $s=m$. Згідно з  \cite[теорема 4.1.4]{KozlovMazyaRossmann97} вимірність коядра цього оператора дорівнює $\dim \mathcal{N}_{\star}$. Оскільки лінійний многовид  $C^{\infty}(\overline{\Omega})\times(C^{\infty}(\Gamma))^{q+\varkappa}$ щільний у просторі $\mathcal{E}^{m-2q,(m+1-2q)}(\Omega,\Gamma)$, то
згідно з \cite[лема 2.1]{HohbergKrein57} існує скінченновимірний простір $$
\mathcal{G}\subset C^{\infty}(\overline{\Omega})\times(C^{\infty}(\Gamma))^{q+\varkappa}
$$
такий, що
\begin{equation}\label{proof--lemma--a}
\mathcal{E}^{m-2q,(m+1-2q)}(\Omega,\Gamma)=\mathcal{G}\dotplus \Lambda(\mathcal{D}^{m,(m+1)}(\Omega,\Gamma)).
\end{equation}
Звідси випливає, що $\dim\mathcal{G}=\dim\mathcal{N}_{\star}$.

Нехай дійсне число $l$ задовольняє умову $m+1/2<l<s$. Тоді виконуються неперервні вкладення
\begin{gather*}
\mathcal{E}^{s-2q,\varphi}(\Omega,\Gamma)\hookrightarrow
\mathcal{E}^{l-2q}(\Omega,\Gamma)=\\
=\mathcal{E}^{l-2q,(m+1-2q)}(\Omega,\Gamma)
\hookrightarrow\mathcal{E}^{m-2q,(m+1-2q)}(\Omega,\Gamma)
\end{gather*}
на підставі \eqref{1f7} і того, що простори $H^{l-2q,(m+1-2q)}(\Omega)$ і $H^{l-2q}(\Omega)$ рівні з точністю до еквівалентності норм при $l-2q>m+1-2q-1/2$, як це зазначалось у доведенні теореми~\ref{1th1}. Окрім того, $\mathcal{G}\subset\mathcal{E}^{s-2q,\varphi}(\Omega,\Gamma)$. Тому з рівності \eqref{proof--lemma--a} випливає формула
\begin{equation}\label{proof--lemma--b}
\mathcal{E}^{s-2q,\varphi}(\Omega,\Gamma)=\mathcal{G}\dotplus
\bigl(\Lambda(\mathcal{D}^{m,(m+1)}(\Omega,\Gamma))\cap
\mathcal{E}^{s-2q,\varphi}(\Omega,\Gamma)\bigr).
\end{equation}
На підставі \cite[теорема 4.1.4]{KozlovMazyaRossmann97} область значень  $\Lambda(\mathcal{D}^{m,(m+1)}(\Omega,\Gamma))$ оператора \eqref{roit--oper}, де $s=m$, складається з усіх векторів
$$
(f,g)\in\mathcal{E}^{m-2q,(m+1-2q)}(\Omega,\Gamma),
$$
які задовольняють умову \eqref{1f11}. Тому другий доданок у сумі \eqref{proof--lemma--b} складається з усіх векторів $(f,g)\in\mathcal{E}^{s-2q,\varphi}(\Omega,\Gamma)$, які задовольняють \eqref{1f11}. Таким чином, формула \eqref{proof--lemma--b} перетворюється на рівність \eqref{1f13}, в якій, згідно з наведеними міркуваннями, простір $\mathcal{G}$ не залежить від $s$ і $\varphi$.

Лему~\ref{1lem1} доведено.

\medskip

\noindent\textbf{Доведення теореми~\ref{1th2}.} Згідно з теоремою~\ref{1th1} звуження оператора \eqref{1f9} на підпростір $\mathcal{P}(\mathcal{D}^{s,\varphi}(\Omega,\Gamma))$ є неперервним і взаємно однозначним відображенням цього підпростору на підпростір $\mathcal{Q}(\mathcal{E}^{s-2q,\varphi}(\Omega,\Gamma))$. Тому на підставі теореми Банаха про обернений оператор це відображення є ізоморфізмом~\eqref{1f14}. Теорему~\ref{1th2} доведено.

\medskip

\noindent\textbf{Доведення теореми~\ref{1th3}.}
У випадку, коли $\chi(\cdot)\equiv\eta(\cdot)\equiv1$ ця теорема є наслідком скінченновимірності ядра і замкненості області значень оператора \eqref{1f9} з теореми~\ref{1th1} та компактності вкладення $\mathcal{D}^{s,\varphi}(\Omega)\hookrightarrow \mathcal{D}^{s-\lambda,\varphi}(\Omega)$ для довільного $\lambda>0$. Це стверджує лема Пітре \cite[лема~3]{Peetre61}. Отже, існує число $\tilde{c}=\tilde{c}(s,\varphi,\lambda)>0$ таке, що для довільного вектора $(u',v')\in \mathcal{D}^{s,\varphi}(\Omega,\Gamma)$ виконується глобальна апріорна оцінка
\begin{equation}\label{global-estimate}
\|(u',v')\|_{\mathcal{D}^{s,\varphi}(\Omega,\Gamma)}\leq\tilde{c}\,
\bigl(\|\Lambda(u', v')\|_{\mathcal{E}^{s-2q,\varphi}(\Omega,\Gamma)}
+\|(u',v')\|_{\mathcal{D}^{s-\lambda,\varphi}(\Omega,\Gamma)}\bigr).
\end{equation}
У цій оцінці число $\lambda>0$ вибране довільним чином.

Виведемо з цієї оцінки теорему~\ref{1th3} для $\lambda=1$. Зауважимо спочатку, що нерівність $\lambda<s-m+1/2$, вказана в умові цієї теореми, виконується для $\lambda=1$. Довільно виберемо вектор $(u,v)\in \mathcal{D}^{s,\varphi}(\Omega,\Gamma)$. Нехай функції $\chi,\eta\in C^{\infty}(\overline{\Omega})$ такі як в умові теореми~\ref{1th3}. Поклавши $(u',v'):=\chi(u,v)\in \mathcal{D}^{s,\varphi}(\Omega,\Gamma)$ та $\lambda:=1$ в оцінці \eqref{global-estimate}, запишемо
\begin{equation}\label{1f29}
\|\chi(u,v)\|_{\mathcal{D}^{s,\varphi}(\Omega,\Gamma)}\leq\tilde{c}\,
\bigl(\|\Lambda(\chi(u,v))\|_{\mathcal{E}^{s-2q,\varphi}(\Omega,
\Gamma)}
+\|\chi(u,v)\|_{\mathcal{D}^{s-1,\varphi}(\Omega,\Gamma)}\bigr).
\end{equation}
Переставивши оператор множення на функцію $\chi$ з
диференціальним оператором $\Lambda$ отримаємо рівність
\begin{equation*}
\Lambda(\chi (u, v))=\Lambda (\chi\eta (u,v))=\chi \Lambda(\eta (u,v))
+\Lambda'(\eta (u,v))=
\chi \Lambda(u,v)+\Lambda'(\eta (u,v)).
\end{equation*}
Тут
$$
\Lambda'(u,v):=\biggl(A'u,\,B'_{1}u+\sum_{k=1}^{\varkappa}C'_{1,k}v_{k},...,
B'_{q+\varkappa}u+\sum_{k=1}^{\varkappa}C'_{q+\varkappa,k}v_{k}\biggr),
$$
де $A'$ --- деякий лінійний диференціальний оператор на $\overline{\Omega}$ порядку $\mathrm{ord}\,A'\leq 2q-1$, $B'_{j}$~--- деякий крайовий лінійний диференціальний оператор на $\Gamma$ порядку $\mathrm{ord}\,B'_{j}\leq m_{j}-1$, а $C'_{j,k}$~--- деякий (дотичний) лінійний диференціальний оператор на $\Gamma$ порядку $\mathrm{ord}\,C'_{j,k}\leq m_{j}+r_{k}-1$. Коефіцієнти цих операторів є нескінченно гладкими функціями на $\overline{\Omega}$ і $\Gamma$ відповідно.

Таким чином,
\begin{equation}\label{1f30}
\Lambda(\chi (u,v))=\chi\Lambda(u,v)+\Lambda'(\eta (u,v)).
\end{equation}
З властивостей порядків компонент диференціального оператора $\Lambda'$ негайно випливає нерівність
\begin{equation}\label{1f30b}
\|\Lambda'(\eta(u,v))\|_{\mathcal{E}^{s-2q,\varphi}(\Omega,\Gamma)}
\leq c_{1}\|\eta(u,v)\|_{\mathcal{D}^{s-1,\varphi}(\Omega,\Gamma)}.
\end{equation}
У цьому доведенні через $c_{1},\ldots,c_{7}$ позначено додатні числа, не залежні від~$(u,v)$.

На підставі формул \eqref{1f29}\,--\,\eqref{1f30b} отримаємо нерівності
\begin{gather*}
\|\chi (u,v)\|_{\mathcal{D}^{s,\varphi}(\Omega,\Gamma)}\leq \tilde{c}\,
\bigl(\|\chi \Lambda(u,v)\|_{\mathcal{E}^{s-2q,\varphi}(\Omega,\Gamma)}+\\
+\|\Lambda'(\eta (u,v))\|_{\mathcal{E}^{s-2q,\varphi}(\Omega,\Gamma)}
+\|\chi(u,v)\|_{\mathcal{D}^{s-1,\varphi}(\Omega,\Gamma)}\bigr)\leq\\
\leq\tilde{c}\,
\|\chi\Lambda(u,v)\|_{\mathcal{E}^{s-2q,\varphi}(\Omega,\Gamma)}+
\tilde{c}\,c_{1}\|\eta (u,v)\|_{\mathcal{D}^{s-1,\varphi}(\Omega,\Gamma)}
+\tilde{c}\,\|\chi (u,v)\|_{\mathcal{D}^{s-1,\varphi}(\Omega,\Gamma)}.
\end{gather*}
Тут
\begin{equation*}
\|\chi (u,v)\|_{\mathcal{D}^{s-1,\varphi}(\Omega,\Gamma)}=
\|\chi\eta (u,v)\|_{\mathcal{D}^{s-1,\varphi}(\Omega,\Gamma)}
\leq c_{2}\|\eta (u,v)\|_{\mathcal{D}^{s-1,\varphi}(\Omega,\Gamma)}.
\end{equation*}
Отже,
\begin{equation}\label{proof-th3}
\|\chi (u,v)\|_{\mathcal{D}^{s,\varphi}(\Omega,\Gamma)}\leq c_{3}
\bigl(\|\chi\Lambda(u,v)\|_{\mathcal{E}^{s-2q,\varphi}(\Omega,\Gamma)}
+\|\eta (u,v)\|_{\mathcal{D}^{s-1,\varphi}(\Omega,\Gamma)}\bigr).
\end{equation}
З цієї нерівності випливає потрібна оцінка \eqref{1f15}, оскільки
\begin{equation*}
\|\chi\Lambda(u,v)\|_{\mathcal{E}^{s-2q,\varphi}(\Omega,\Gamma)}=
\|\chi\eta\Lambda(u,v)\|_{\mathcal{E}^{s-2q,\varphi}(\Omega,\Gamma)}
\leq c_{4}\|\eta\Lambda(u,v)\|_{\mathcal{E}^{s-2q,\varphi}(\Omega,\Gamma)}.
\end{equation*}
Теорему \ref{1th3} доведено у випадку, коли $\lambda=1$. Звісно, її висновок правильний і якщо $0<\lambda<1$.

Доведемо тепер цю теорему у випадку, коли
\begin{equation}\label{proof-th3-a}
1<\lambda<s-m+1/2.
\end{equation}
Для кожного дійсного числа $l\geq1$ позначимо через $\mathcal{K}_{l}$ твердження теореми~\ref{1th3} у випадку, коли $\lambda=l$ і фіксоване $\varphi\in\mathcal{M}$. А саме, $\mathcal{K}_{l}$ позначає таке твердження: для довільних числа $s>m+1/2$ і функцій $\chi,\eta\in C^{\infty}(\overline{\Omega})$, які задовольняють умови $l<s-m+1/2$ і $\eta=1$ в околі $\mathrm{supp}\,\chi$, існує число $c=c(s,\varphi,l,\chi,\eta)>0$ таке, що для довільного вектора $(u,v)\in \mathcal{D}^{s,\varphi}(\Omega,\Gamma)$ виконується нерівність \eqref{1f15} з $\lambda=l$. Істинність твердження $\mathcal{K}_{1}$ доведена вище. Довільно виберемо дійсні числа $l\geq1$ і $\delta\in(0,1]$. Доведемо, що $\mathcal{K}_{l}\Rightarrow\mathcal{K}_{l+\delta}$.

Припустимо, що твердження $\mathcal{K}_{l}$ істинне. Нехай число $s>m+1/2$ і функції $\chi,\eta\in C^{\infty}(\overline{\Omega})$ задовольняють умови $l+\delta<s-m+1/2$ і $\eta=1$ в околі $\mathrm{supp}\,\chi$. Тоді знайдеться функція $\eta_{1}\in C^{\infty}(\overline{\Omega})$ така, що $\eta_{1}=1$ в околі $\mathrm{supp}\,\chi$ і $\eta=1$ в околі $\mathrm{supp}\,\eta_{1}$. За припущенням, існує число $c_{5}>0$ таке, що для довільного вектора $(u,v)\in \mathcal{D}^{s,\varphi}(\Omega,\Gamma)$ виконується оцінка
\begin{equation}\label{proof-th3-b}
\|\chi (u,v)\|_{\mathcal{D}^{s,\varphi}(\Omega,\Gamma)}\leq c_{5}
\bigl(\|\eta_{1}\Lambda(u,v)\|_{\mathcal{E}^{s-2q,\varphi}(\Omega,\Gamma)}
+\|\eta_{1} (u,v)\|_{\mathcal{D}^{s-l,\varphi}(\Omega,\Gamma)}\bigr).
\end{equation}
Оскільки $s-l-\delta+1>m+1/2$, то на підставі твердження $\mathcal{K}_{1}$ маємо оцінку
\begin{equation}\label{proof-th3-c}
\begin{gathered}
\|\eta_{1}(u,v)\|_{\mathcal{D}^{s-l,\varphi}(\Omega,\Gamma)}\leq
\|\eta_{1}(u,v)\|_{\mathcal{D}^{s-l-\delta+1,\varphi}(\Omega,\Gamma)}\leq\\
\leq c_{6}\bigl(\|\eta\Lambda(u,v)\|_{\mathcal{E}^{s-l-\delta+1-2q,\varphi}(\Omega,\Gamma)}
+\|\eta (u,v)\|_{\mathcal{D}^{s-l-\delta,\varphi}(\Omega,\Gamma)}\bigr).
\end{gathered}
\end{equation}
Окрім того,
\begin{equation}\label{proof-th3-d}
\|\eta_{1}\Lambda(u,v)\|_{\mathcal{E}^{s-2q,\varphi}(\Omega,\Gamma)}=
\|\eta_{1}\eta\Lambda(u,v)\|_{\mathcal{E}^{s-2q,\varphi}(\Omega,\Gamma)}\leq c_{7}\|\eta\Lambda(u,v)\|_{\mathcal{E}^{s-2q,\varphi}(\Omega,\Gamma)}.
\end{equation}
На підставі оцінок \eqref{proof-th3-b}\,--\,\eqref{proof-th3-d} запишемо
\begin{gather*}
\|\chi (u,v)\|_{\mathcal{D}^{s,\varphi}(\Omega,\Gamma)}\leq c_{5}c_{7}
\|\eta\Lambda(u,v)\|_{\mathcal{E}^{s-2q,\varphi}(\Omega,\Gamma)}+\\
+c_{5}c_{6}\bigl(\|\eta\Lambda(u,v)\|_
{\mathcal{E}^{s-2q,\varphi}(\Omega,\Gamma)}
+\|\eta(u,v)\|_{\mathcal{D}^{s-l-\delta,\varphi}(\Omega,\Gamma)}\bigr),
\end{gather*}
тобто отримано нерівність \eqref{1f15} з $\lambda=l+\delta$. Імплікація  $\mathcal{K}_{l}\Rightarrow\mathcal{K}_{l+\delta}$ обґрунтована.

Тепер можемо довести теорему~\ref{1th3} у випадку \eqref{proof-th3-a}. За доведеним, правильний ланцюжок імплікацій
\begin{equation*}
\mathcal{K}_{1}\Rightarrow\mathcal{K}_{2}\Rightarrow\ldots
\Rightarrow\mathcal{K}_{[\lambda]}\Rightarrow \mathcal{K}_{\lambda},
\end{equation*}
де твердження $\mathcal{K}_{1}$ істинне, а $\mathcal{K}_{\lambda}$ є твердженням теореми~\ref{1th3} у досліджуваному випадку (як звичайно, $[\lambda]$~--- ціла частина числа~$\lambda$).

Теорему~\ref{1th3} доведено.

У зауваженні~\ref{rem1} потребують обґрунтування друге і останнє речення. Друге речення обґрунтоване у першому абзаці доведення цієї теореми, а останнє речення є прямим наслідком оцінки \eqref{proof-th3}.

\medskip

\noindent\textbf{Доведення теореми~\ref{1th4}.} Спочатку доведемо цю теорему у глобальному випадку, коли $\Omega_{0}=\Omega$ і $\Gamma_{0}=\Gamma$. За умовою,
$(u,v)\in \mathcal{D}^{l}(\Omega,\Gamma)$ для деякого дійсного числа $l>m+1/2$ і, окрім того, $(f,g)\in\mathcal{E}^{s-2q,\varphi}(\Omega,\Gamma)$. Тому
$$
(f,g)=\Lambda(u,v)\in\mathcal{E}^{s-2q,\varphi}(\Omega,\Gamma)
\cap\Lambda(\mathcal{D}^{l}(\Omega,\Gamma))=
\Lambda(\mathcal{D}^{s,\varphi}(\Omega,\Gamma))
$$
на підставі теореми~\ref{1th1}. Отже, поряд з умовою $\Lambda(u,v)=(f,g)$ виконується рівність $\Lambda(u',v')=(f,g)$ для деякого
$(u',v')\in\mathcal{D}^{s,\varphi}(\Omega,\Gamma)$. Тому $\Lambda(u-u',v-v')=0$, де $(u-u',v-v')\in\mathcal{D}^{l}(\Omega,\Gamma)$. Це за теоремою~\ref{1th1}, у якій беремо $l$ замість $s$ і покладаємо $\varphi(\cdot)\equiv1$, тягне за собою включення
$$
(u-u',v-v')\in \mathcal{N}\subset C^{\infty}(\overline{\Omega})\times (C^{\infty}(\Gamma))^{\varkappa}.
$$
Отже,
$$
(u,v)=(u',v')+(u-u',v-v')\in\mathcal{D}^{s,\varphi}(\Omega,\Gamma).
$$
Теорему~\ref{1th4} доведено у досліджуваному випадку.

Виведемо з цього випадку теорему~\ref{1th4} у загальній ситуації. Попередньо доведемо, що за умови цієї теореми виконується така імплікація для кожного дійсного числа $\sigma>m+1/2$:
\begin{equation}\label{proof4-impl}
(u,v)\in\mathcal{D}^{\sigma,\varphi}_{\mathrm{loc}}(\Omega_{0},\Gamma_{0})
\,\Rightarrow\,
(u,v)\in
\mathcal{D}^{\min\{\sigma+1,s\},\varphi}_
{\mathrm{loc}}(\Omega_{0},\Gamma_{0}).
\end{equation}

Припустимо, що виконується посилка цієї імплікації для деякого $\sigma>m+1/2$. Виберемо довільним чином функції $\chi$ і $\eta$ такі як у теоремі~\ref{1th3}. Згідно з формулою \eqref{1f30} маємо рівність
\begin{equation*}
\Lambda(\chi (u,v))=\chi(f,g)+\Lambda'(\eta(u,v)).
\end{equation*}
Тут $\chi(f,g)\in\mathcal{E}^{s-2q,\varphi}(\Omega,\Gamma)$ за умовою і $$
\Lambda'(\eta(u,v))\in\mathcal{E}^{\sigma-2q+1,\varphi}(\Omega,\Gamma)
$$
на підставі посилки імплікації та властивостей порядків компонент оператора $\Lambda'$, зазначених у доведенні теореми~\ref{1th3}.
Отож,
$$
\Lambda(\chi (u,v))\in\mathcal{E}^{\min\{\sigma+1,s\}-2q,\varphi}(\Omega,\Gamma).
$$
Звідси, за теоремою~\ref{1th4} у вже доведеному глобальному випадку, випливає включення
$$
\chi (u,v)\in\mathcal{D}^{\min\{\sigma+1,s\},\varphi}(\Omega,\Gamma).
$$
Таким чином,
$$
(u,v)\in\mathcal{D}^{\min\{\sigma+1,s\},\varphi}_
{\mathrm{loc}}(\Omega_{0},\Gamma_{0})
$$
з огляду на довільність нашого вибору функції $\chi$. Імплікацію \eqref{proof4-impl} доведено.

Оскільки $(u,v)\in\mathcal{D}^{m+1/2+}(\Omega,\Gamma)$ за умовою, то
$$
(u,v)\in\mathcal{D}^{l,\varphi}(\Omega,\Gamma)\subset
\mathcal{D}^{l,\varphi}_{\mathrm{loc}}(\Omega_{0},\Gamma_{0})
$$
для деякого числа $l\in(m+1/2,s)$. Застосувавши імплікацію \eqref{proof4-impl} $[s-l]+1$ разів послідовно для значень $\sigma=l$, $\sigma=l+1$, ..., $\sigma=l+[s-l]$, отримаємо такі співвідношення:
\begin{gather*}
(u,v)\in\mathcal{D}^{l,\varphi}_{\mathrm{loc}}(\Omega_{0},\Gamma_{0})\,
\Rightarrow\\
\Rightarrow\,(u,v)\in\mathcal{D}^{\min\{l+1,s\},\varphi}_
{\mathrm{loc}}(\Omega_{0},\Gamma_{0})=
\mathcal{D}^{l+1,\varphi}_{\mathrm{loc}}(\Omega_{0},\Gamma_{0})\,
\Rightarrow\ldots\\ \Rightarrow\,
(u,v)\in\mathcal{D}^{\min\{l+[s-l]+1,s\},\varphi}_
{\mathrm{loc}}(\Omega_{0},\Gamma_{0})=
\mathcal{D}^{s,\varphi}_{\mathrm{loc}}(\Omega_{0},\Gamma_{0}).
\end{gather*}

Теорему~\ref{1th4} доведено.

\medskip

\noindent\textbf{Доведення теореми~\ref{1th5}.} Розглянемо спочатку випадок, коли $l+\nobreak n/2>m+1/2$. Тоді, на підставі теореми~\ref{1th4}, з умови  \eqref{1th5-cond} випливає включення $u\in H^{l+n/2,\varphi}_{\mathrm{loc}}(\Omega_{0},\Gamma_{0})$. Довільно виберемо точку
$x\in\Omega_{0}\cup\Gamma_{0}$ і функцію $\chi\in C^{\infty}(\overline{\Omega})$ таку, що $\mathrm{supp}\,\chi\subset\Omega_0\cup\Gamma_0$ і $\chi=1$ у деякому околі $V(x)$ (у топології на $\overline{\Omega}$) точки~$x$. Згідно з твердженням~\ref{1prop1} маємо включення $\chi u\in H^{l+n/2,\varphi}(\Omega)\subset C^{l}(\overline{\Omega})$.
Тому $u\in C^{l}(V(x))$. Звідси, з урахуванням довільності вибору точки $x\in\Omega_{0}\cup\Gamma_{0}$, робимо висновок, що $u\in C^{l}(\Omega_{0}\cup\Gamma_{0})$.

Дослідимо тепер випадок, коли $l+\nobreak n/2\leq m+1/2$. Оскільки, за умовою, $(u,v)\in\mathcal{D}^{m+1/2+}(\Omega,\Gamma)$, то
$u\in H^{m+1/2+\varepsilon}(\Omega)$ для деякого $\varepsilon>0$. Отож, у цьому випадку $u\in H^{l+n/2+\varepsilon}(\Gamma)\subset C^{l}(\overline{\Omega})$ за теоремою вкладення Соболєва. Отже, і поготів $u\in C^{l}(\Omega_{0}\cup\Gamma_{0})$.

Теорему~\ref{1th5} доведено.

\medskip

\noindent\textbf{Доведення теореми~\ref{1th6}.} Розглянемо спочатку випадок, коли $l-\nobreak r_{k}+n/2>m+1/2$. Тоді, на підставі теореми~\ref{1th4}, з умови \eqref{1th6-cond} випливає включення $(u,v)\in\mathcal{D}^{l-r_{k}+n/2,\varphi}_
{\mathrm{loc}}(\Omega_{0},\Gamma_{0})$. Зокрема, $v_{k}\in H^{l+(n-1)/2,\varphi}_{\mathrm{loc}}(\Gamma_{0})$.
Довільно виберемо точку $x\in\Gamma_{0}$ і функцію $\chi\in C^{\infty}(\Gamma)$ таку, що $\mathrm{supp}\,\chi\subset\Gamma_{0}$ і $\chi=1$ в деякому околі $V(x)\subset\Gamma$ точки $x$. На підставі твердження~\ref{1prop1} маємо включення
$\chi v_{k}\in H^{l+(n-1)/2,\varphi}(\Gamma)\subset C^{l}(\Gamma)$.
Отже, $v_{k}\in C^{l}(V(x))$. Звідси, з урахуванням довільності вибору точки $x\in\Gamma_{0}$, випливає потрібна властивість $v_{k}\in C^{l}(\Gamma_{0})$.

Дослідимо тепер випадок, коли $l-r_{k}+n/2\leq m+1/2$. Оскільки, за умовою, $(u,v)\in\mathcal{D}^{m+1/2+}(\Omega,\Gamma)$, то $v_{k}\in H^{m+r_{k}+\varepsilon}(\Omega)$ для деякого $\varepsilon>0$. Отже, у цьому випадку $v_{k}\in H^{l+(n-1)/2+\varepsilon}(\Gamma)\subset C^{l}(\Gamma)$ за теоремою вкладення Соболєва. Отож, і поготів $v_{k}\in C^{l}(\Gamma_{0})$.

Теорему~\ref{1th6} доведено.

\medskip

\noindent\textbf{Доведення теореми~\ref{1th7}.} Включення $u\in C^{2q}(\Omega)$ є наслідком умови \eqref{1f46} на підставі теореми 5, у якій покладаємо $l:=2q$, $\varphi:=\nobreak\varphi_1$, $\Omega_{0}:=\Omega$ і $\Gamma_{0}:=\varnothing$. Включення $u\in C^{m}(U_{\sigma}\cup\Gamma)$ є наслідком умови \eqref{1f46} на підставі теореми 5, у якій беремо $l:=m$, $\varphi:=\varphi_1$, $\Omega_{0}:=U_{\sigma}$ i $\Gamma_{0}:=\Gamma$. Включення $v_{k}\in C^{m+r_k}(\Gamma)$ є наслідком умови \eqref{1f47} на підставі теореми 6, у якій беремо $l:=m+r_k$, $\varphi:=\varphi_2$ i $\Gamma_{0}:=\Gamma$. Таким чином, розв'язок $(u,v)$ класичний.

Теорему 7 доведено.

На завершення цього пункту обґрунтуємо зауваження~\ref{rem2}. Нехай $0\leq l\in\mathbb{Z}$ і $\varphi\in\mathcal{M}$. Доведемо  спочатку еквівалентність $\eqref{1f45}\Leftrightarrow\eqref{rem1-impl-a}$ за умови $l+n/2>m+1/2$. З огляду на теорему \ref{1th5} залишається показати, що $\eqref{rem1-impl-a}\Rightarrow\eqref{1f45}$. Припустимо, що $\eqref{rem1-impl-a}$ істинне. Нехай $\Omega_{1}$~--- деяка відкрита куля у просторі $\mathbb{R}^{n}$, замикання якої лежить в $\Omega_{0}$. Для довільного розподілу $w\in H^{l+n/2,\varphi}(\mathbb{R}^{n})$  розглянемо вектор $(u,v)$, де $u:=w\!\upharpoonright\!\Omega$, а $v:=0$. Вектор
$$
(u,v)\in\mathcal{D}^{l+n/2,\varphi}(\Omega,\Gamma)\subset
\mathcal{D}^{m+1/2+}(\Omega,\Gamma)
$$
задовольняє умову теореми~\ref{1th5}, а тому $u\in C^{l}(\Omega_{0}\cup\Gamma_{0})$ згідно із зробленим припущенням. Отже,
$$
H^{l+n/2,\varphi}(\Omega_{1})=\bigl\{w\!\upharpoonright\!\Omega_{1}:
w\in H^{l+n/2,\varphi}(\mathbb{R}^{n})\bigr\}\subset C^{l}(\overline{\Omega_{1}}).
$$
Це вкладення тягне за собою умову \eqref{1f45} на підставі твердження~\ref{1prop1}, де замість $\Omega$ беремо кулю $\Omega_{1}$. Таким чином, $\eqref{rem1-impl-a}\Rightarrow\eqref{1f45}$.

Припустимо тепер, що $k\in\{1,...,\varkappa\}$,  $\Gamma_{0}\neq\varnothing$ і $l-r_{k}+n/2>m+1/2$. Доведемо  еквівалентність $\eqref{1f45}\Leftrightarrow\eqref{rem1-impl-b}$. З огляду на теорему \ref{1th6} залишається показати, що $\eqref{rem1-impl-b}\Rightarrow\eqref{1f45}$. Припустимо, що $\eqref{rem1-impl-b}$ істинне. Нехай $W$~--- деяка відкрита непорожня підмножина межі $\Gamma$ така, що $\overline{W}\subset\Gamma_{0}\cap\Gamma_{j}$ для деякого $j\in\{1,\ldots,\lambda\}$. Тут, нагадаємо, $\pi_{j}:\mathbb{R}^{n-1}\leftrightarrow\Gamma_{j}$ є локальною картою на многовиді $\Gamma$, яка використана в означенні просторів Хермандера на~$\Gamma$. Розглянемо довільний розподіл $h\in H^{l+(n-1)/2,\varphi}(\mathbb{R}^{n-1})$ такий, що $\mathrm{supp}\,h\subset U$, де $U:=\pi_{j}^{-1}(W)$~--- відкрита підмножина простору $\mathbb{R}^{n-1}$. Розглянемо вектор $(u,v):=(u,v_{1},\ldots,v_{\varkappa})$, де $u:=0$, $v_{k}:=\Theta_{j}(h\circ\pi_{j}^{-1})$ і $v_{i}:=0$ при $i\neq k$. Тут $\Theta_{j}$~--- оператор продовження нулем фінітного розподілу з підмножини $\Gamma_{j}$ на увесь многовид $\Gamma$. Вектор
$$
(u,v)\in\mathcal{D}^{l-r_{k}+n/2,\varphi}(\Omega,\Gamma)
\subset\mathcal{D}^{m+1/2+}(\Omega,\Gamma)
$$
задовольняє умову теореми~\ref{1th6} й тому $v_{k}\in C^{l}(\Gamma_{0})$ згідно зробленого припущення. Отже, $h=v_{k}\circ\pi_{j}\in C^{l}(\mathbb{R}^{n-1})$. Нехай $U_{1}$~--- відкрита підмножина простору $\mathbb{R}^{n-1}$ така, що $\overline{U_{1}}\subset U$. Оскільки простір $H^{l+(n-1)/2,\varphi}(U_{1})$ складається зі звужень на $U_{1}$ усіх розглянутих розподілів $h$, то
$$
H^{l+(n-1)/2,\varphi}(U_{1})\subset C^{l}(\overline{\Omega_{1}}).
$$
Це вкладення тягне за собою умову \eqref{1f45} на підставі твердження~\ref{1prop1}, де замість $\Omega$ беремо кулю $U_{1}$, а замість $n$ беремо $n-1$. Таким чином, $\eqref{rem1-impl-b}\Rightarrow\eqref{1f45}$.

Зауваження~\ref{rem2} обґрунтоване.

\section{Висновки}\label{1sec7}

У статті досліджено еліптичну за Лавруком задачу з крайовими операторами високих порядків в уточненій соболєвській шкалі. Доведено, що цій задачі відповідають нетерові обмежені оператори, які діють у відповідних парах просторів Хермандера (теорема~\ref{1th1}) та породжують ізоморфізми між їх підпросторами скінченної ковимірності (теорема~\ref{1th2}). Встановлено локальну апріорну оцінку узагальнених розв'язків задачі (теорема~\ref{1th3}). Досліджено регулярність узагальнених розв'язків в уточненій соболєвській шкалі (теорема~\ref{1th4}). Як застосування уточненої соболєвської шкали, знайдено нові достатні умови неперервності узагальнених похідних (заданого порядку) розв'язків досліджуваної задачі (теореми \ref{1th5} і \ref{1th6}) та умову класичності узагальненого розв'язку задачі (теорема~\ref{1th7}).

\medskip

\textit{Автори дякують О. О. Мурачу за керівництво роботою.}

\medskip

\end{document}